\documentclass[11pt,reqno]{amsart}

\usepackage[a4paper,margin=2.4cm]{geometry}

\usepackage{amssymb}
\usepackage{amsmath}
\usepackage{amsthm}

\usepackage{url}
\usepackage{xcolor} 
\usepackage{latexsym}
\usepackage{amsfonts}
\usepackage{mathtools}
\usepackage{enumerate}
\usepackage{mathabx}
\usepackage{tikz}
\usepackage{graphicx}
\usepackage{microtype}
\usepackage{placeins}
\usepackage{float}
\usepackage[calcwidth=0.8\linewidth]{caption}
\usepackage{algorithm}
\usepackage{algpseudocode}
\usetikzlibrary {arrows.meta} 
\usepackage{pgfplots}
\usetikzlibrary{shapes.geometric, arrows}
\usepackage{nicefrac}

\newif\ifrevision
\revisionfalse
\definecolor{revcolor}{rgb}{0.75,0.0,0.0}
\ifrevision
\newcommand{\new}[1]{{\color{revcolor}#1}} %
\else
\newcommand{\new}[1]{#1}
\fi

\tikzstyle{box} = [rectangle, rounded corners, minimum width=4cm, minimum height=.7cm,text centered, draw=black, fill=black!30]
\tikzstyle{midbox} = [rectangle, rounded corners, minimum width=2.62cm, minimum height=.7cm,text centered, draw=black, fill=black!30]
\tikzstyle{arrow} = [thick,->,>=stealth]

\algnewcommand\algorithmicgoal{\textbf{Goal:}}
\algnewcommand\Goal{\item[\algorithmicgoal]}
\algnewcommand\algorithmicinput{\textbf{Input:}}
\algnewcommand\Input{\item[\algorithmicinput]}
\algnewcommand\algorithmicoutput{\textbf{Output:}}
\algnewcommand\Output{\item[\algorithmicoutput]}
\algnewcommand\algorithmicbreak{\textbf{Break:}}
\algnewcommand\Break{\item[\algorithmicbreak]}

\newtheorem{theorem}{Theorem}[section]
\newtheorem{remark}[theorem]{Remark}

\newtheorem{assumption}[theorem]{Assumption}
\newtheorem{proposition}[theorem]{Proposition}
\newtheorem{definition}[theorem]{Definition}
\newtheorem{lemma}[theorem]{Lemma}
\newtheorem{corollary}[theorem]{Corollary}

\newcommand{\R}{\ensuremath{\mathbb R}}    %

\providecommand{\ddt}{\frac{\mathrm{d}}{\mathrm{d}t}}
\providecommand{\at}[2]{\ensuremath{\left. #1 \right\vert_{#2}}}

\newcommand{\Hess}{\operatorname{Hess}}
\newcommand{\grad}{\operatorname{grad}}

\newcommand{\im}{\operatorname{im}}
\newcommand{\Flow}{\operatorname{Fl}}
\newcommand{\argmin}{\operatorname{arg} \hspace*{-0.5cm}\min}
\newcommand{\one}{\mathbf{1}}
\newcommand{\tD}{\mathrm{D}}
\newcommand{\tT}{\mathrm{T}}

\newcommand{\karl}[1]{#1}
\newcommand{\kw}[1]{}

\begin{document}
	
	\title[Maximum-likelihood reprojections for reliable Koopman-based prediction]{\new{Maximum-likelihood reprojections for reliable Koopman predictors of parametric nonlinear %
			dynamical systems}}
	
	\author{Pieter van Goor}
	\address{School of Aerospace, Mechanical, and Mechatronic Engineering, The University of Sydney, NSW 2006, Australia}
	
	\author{Manuel Schaller}
	\address{Faculty of Mathematics, Technische Universit\"at Chemnitz, 09126 Chemnitz, Germany}
	\email{manuel.schaller@math.tu-chemnitz.de}
	
	\author{Robert Mahony}
	\address{Systems Theory and Robotics Group, School of Engineering, The Australian National University, Canberra 0200, Australia}
	
	\author{Karl Worthmann}
	\address{Optimization-based Control Group, Institute of Mathematics, Technische Universit\"at Ilmenau, 98693 Ilmenau, \new{Germany}}
	
	\thanks{K.\ Worthmann gratefully acknowledges funding by the Carl-Zeiss foundation within the project KI-MSO-O}

	
	
	\begin{abstract}
		Koopman-based methods leverage a nonlinear lifting to enable linear regression techniques. Consequently, \karl{system analysis}, %
		learning, and prediction are \karl{considered} %
		through the lens of this lifting \karl{resulting in} %
		a nonlinear manifold that is invariant under the \karl{action of the} Koopman operator. 
		\karl{Typically, this invariance is not preserved in d}ata-driven approximations such as Extended Dynamic Mode Decomposition (EDMD). %
		In this work, \karl{we systematically generalize the coordinate reprojection to render Koopman predictors more reliable and versatile. 
			To this end, we propose consistency projectors to align the dynamics of the Koopman surrogate model with the geometry of the nonlinear manifold using %
			concepts from differential geometry. Further, we establish a link to %
			maximum likelihood estimates, which can be leveraged as an event-triggering mechanism to only reproject if really necessary. 
			In addition, we extend %
			our approach %
			to parameter-affine systems enabling sensitivity analysis}. 
		In numerical examples ranging from a supercritical pitchfork bifurcation over the Duffing oscillator to the Lorenz system, we showcase that the reprojected \karl{Koopman} surrogate %
		reproduces the bifurcation behaviour of the underlying system, i.e., the emerging equilibria, their stability and the resulting attractors, whereas predictions without reprojection either converge to spurious equilibria or diverge. This provides the basis for data-driven bifurcation analysis and control applications.
	\end{abstract}
	
	\maketitle

	\section{Introduction}\label{sec:intro}
	\noindent Lifting systems to higher dimensions is a key component of various data-driven methods. In Koopman-based approaches such as Extended Dynamic Mode Decomposition (EDMD), this lifting is performed by means of a dictionary of nonlinear functions~\cite{WillMatt15}, see also the recent overview of EDMD variants~\cite{Colb24}. EDMD serves as a powerful data-driven tool for prediction~\cite{BrunBudi21,GianDas20}, control~\cite{KordMezi18b,BevaSosn2021,MaurMezi20,BoldGrun23,strasser2024koopman,StraWort26}, model reduction \cite{KlusNusk20} or spectral analysis \cite{Mezi05,RowlMezi09,Mezi13,Mezi20,colbrook2024rigorous,FroyGott16} revealing key %
	properties and patterns (such as attractors) of the underlying nonlinear system, see also the recent tutorial paper~\cite{mezic2026koopman}.
	
	In Koopman-based approaches, the \karl{finite dictionary used to lift the nonlinear %
		system %
		yields} %
	a nonlinear manifold, i.e., the injective 
	image of the original finite-dimensional state-space. 
	This manifold is invariant under the action of the Koopman operator. %
	\karl{Hence, %
		analysis, learning, and prediction are conducted solely %
		through the lens of the lifting functions contained in the dictionary}. 
	However, due to errors in data-driven and finite-dimensional approximations, this invariance is usually not preserved, \karl{which can
		lead to %
		instabilities and errors accumulation in predictions of quantities of interest along the flow of the nonlinear dynamical system. Consequently}, 
	in numerical implementations, a (re-)projection step is \karl{often implicitly incorporated} %
	to maintain consistency, specifically by projecting back onto the manifold~\cite{MaurGonc16,MaurGonc19}, thereby mitigating deteriorating effects such as projection errors stemming from finitely many observable functions. %
	After each prediction step of the linear model in the lifted space, the propagated observables are projected back onto the manifold before the surrogate model is used for the next prediction. %
	When the coordinate functions are included in the dictionary, there is a canonical choice for this projection given by reading off the original state-coordinates from the corresponding Koopman coordinate functions~\cite{KordMezi18b,BrunBrun16,PeitOtto20}. 
	Although this reprojection is conceptually simple, it has two limitations. First, a state-inclusive dictionary may be difficult to realize, both from a practical point of view, e.g., if the lifting functions correspond to sensor measurements that do not comprise the full state, and from a conceptual point of view, e.g., for time-delay or distributed parameter systems, whose state is infinite-dimensional. Second, more sophisticated reprojections may yield more accurate predictions and thus allow to reduce the number of reprojection steps required to maintain a desired approximation accuracy. However, to the best of the authors' knowledge, a systematic analysis of reprojection methods is not yet available.
	
	In this work, we \karl{propose a systematic framework for} %
	reprojections onto the \karl{Koopman-induced} manifold \karl{using concepts from differential geometry}. %
	More precisely, we consider distances induced by weighting matrices that correspond to Riemannian metrics on the manifold. 
	The widely-used coordinate projection \karl{is recovered as a special case and} can be interpreted as a closest-point reprojection with respect to a singular weighting matrix.
	\karl{Further}, %
	we propose a weighting matrix that estimates the covariance of the underlying prediction error, allowing an interpretation as maximum likelihood estimates. 
	\karl{This allows to infer event-triggering mechanisms to apply the consistency projector (resulting in a \textit{nonlinear} reprojection step) only if really necessary.}. \karl{Furthermore, we extend the proposed consistency projectors %
		to enable a sensitivity analysis tailored to parameter-affine systems. Finally}, 
	we address the reliable \emph{prediction} of bifurcation behaviour, i.e., for a set of fixed parameter values, we investigate whether the surrogate reproduces the corresponding %
	equilibria, their stability behaviour, and the respective %
	attractors. Detecting the bifurcation point from data, parameter continuation, and a spectral analysis \karl{based on the derived parameter-affine Koopman} %
	surrogate %
	are left to future work. %

\medskip
The contribution of this paper is a geometric framework for reprojections for \karl{parameter-affine} %
systems addressing the following:
\begin{itemize}
	\item[(i)] Analysis and interpretation of the geometric projection in view of maximum likelihood reprojections (Subsection~\ref{subsec:reprojections})\footnote{\karl{This was only briefly sketched in our prior conference paper~\cite{vanGMaho23}.}}.
	\item [(ii)] \karl{Performance measure leveraging a recursion formula for the co-variance propagation, which might be used as an event-triggering mechanism to reduce the number of reprojection steps in multi-step Koopman predictions.}
	\item [(iii)] Extensive numerical case study in view of attractors and bifurcations in parametric systems (Subsection~\ref{subsec:numres}).
	\item[(iv)] Efficient implementation via adaptive reprojections in multi-step predictions (Subsection~\ref{subsec:multistep}) and a Riemannian Newton method for the reprojection step (Subsection~\ref{subsec:newton}).
\end{itemize}

\medskip
\new{\textbf{Relation to the prior conference paper~\cite{vanGMaho23}.} This work builds upon our conference paper~\cite{vanGMaho23}, where reprojections were introduced for \emph{autonomous} systems. Therein, we proposed the weighted closest-point reprojection of Definition~\ref{def:repr} and the interpretation of the widely-used coordinate projection as the special case of a singular weighting matrix, cf.\ Remark~\ref{rem:coord}. Both are recalled in Section~\ref{sec:problem_formulation}. Moreover, the maximum-likelihood interpretation for autonomous systems was only sketched in~\cite{vanGMaho23} and is rigorously stated and proven in Proposition~\ref{prop:maxlikely}. All remaining results are novel, i.e., the parameter-affine (bilinear) Koopman representation~\eqref{eq:parameter_koopman}--\eqref{eq:parameter_kooopman_model}, the parameter-affine covariance surrogate~$Q$ of Proposition~\ref{prop:cov} together with the parametric maximum-likelihood reprojection of Corollary~\ref{cor:mlparam}, the adaptive multistep scheme of Subsection~\ref{subsec:multistep}, the Riemannian Newton method of Subsection~\ref{subsec:newton} as well as the numerical study of Subsection~\ref{subsec:numres}.}

\medskip
\textbf{Outline.}  In Section~\ref{sec:edmd} we provide the fundamentals of Koopman-based approximations. 
In Section~\ref{sec:problem_formulation}, we recall reprojections onto the image of the dictionary to ensure consistency in data-driven approximations as outlined in \cite{vanGMaho23}. 
The extension to the case of parameter-affine systems is provided in Section~\ref{sec:parameter-affine}, which contains the main novel theoretical contributions of this work as well as a numerical study. Last, in Section~\ref{sec:efficient}, we provide an efficient implementation of the proposed method leveraging multistep predictions and a Riemannian Newton method.

\section{The Koopman operator of a dynamical system and its compression}
\label{sec:edmd}

\noindent For an autonomous, locally-Lipschitz continuous vector field $f : \mathbb{R}^d \rightarrow \mathbb{R}^d$, the flow of~$f$ 
is denoted by $\Flow_f^t : \mathbb{R}^d \to \mathbb{R}^d$ on its maximal interval of existence~$I$, i.e., for every $x \in \mathbb{R}^d$, the flow map is defined by
\begin{equation}\nonumber
	\Flow_f^0(x) = x, \qquad
	\ddt \Flow_f^t(x) = f(\Flow_f^t(x)) \qquad\forall\,t \in I.
\end{equation}
Furthermore, let the compact and convex set $\mathbb{X} \subset \mathbb{R}^d$ be forward invariant with respect to the flow of~$f$, i.e., the inclusion $\Flow_f^t(x) \in \mathbb{X}$ holds for all $t \in I = [0,\infty)$.
Then, the strongly-continuous Koopman semigroup of linear operators $(\mathcal{K}_f^t)_{t\geq 0}$ on $L^2(\mathbb{X},\mathbb{R})$ is defined via the identity
\begin{equation}\label{eq:Koopman}
	(\mathcal{K}_f^t \varphi)(x) = \varphi(\Flow_f^t(x)) \qquad\forall\,x \in \mathbb{X}, \varphi \in L^2(\mathbb{X},\mathbb{R})
\end{equation}
for all $t \in [0,\infty)$.
We refer to elements $\varphi \in L^2(\mathbb{X},\mathbb{R})$ as \emph{observables}. 
The forward-invariance assumption on the set~$\mathbb{X}$ is imposed (and in the following tacitly used) to keep the presentation technically simpler, see, e.g., \cite{ZhanZuaz23} for a detailed discussion. 
Otherwise, the Koopman operator~$\mathcal{K}_f^t$ maps observables defined on $f(\mathbb{X}) := \{ y \in \mathbb{R}^d \mid \exists\,x \in \mathbb{X}: \Flow_f^t(x) = y \}$, i.e., the image of the set~$\mathbb{X}$ under the flow, to observables defined on the set~$\mathbb{X}$, see, e.g., \cite{KohnPhil24} for details and ramifications with respect to multi-step predictions.

We briefly introduce EDMD as a method to compute a finite-dimensional approximation of the Koopman operator~$\mathcal{K} = \mathcal{K}_f^t$ for an arbitrary, but fixed time step $t \in (0,\infty)$ and a given vector field~$f$. 
To this end, one selects~$M$ linearly-independent observables $\psi_k \in L^2(\mathbb{X},\mathbb{R})$, $k \in [1:M]$.
Then, we refer to $\mathbb{V} := \operatorname{span} \{ \psi_i \mid i \in [1:M] \}$ as the \emph{dictionary} and define the vector-valued function
\begin{align}\label{eq:Psi}
	\Psi:\mathbb{X} \to \mathbb{R}^M, \qquad x \mapsto (\psi_1(x)\,\ldots\,\psi_M(x))^\top
\end{align}
such that for any $\varphi \in \mathbb{V}$, there is $a\in \R^M$ with $\varphi(x) = a^\top \Psi(x)$.
The application of the Koopman operator $\mathcal{K}$ to $\Psi$ is to be understood componentwise, i.e., we set
\begin{align*}
	\mathcal{K} \Psi = \begin{pmatrix}
		\mathcal{K} \psi_1\\
		\vdots\\
		\mathcal{K} \psi_M
	\end{pmatrix} = \Psi(\Flow_f^t(\cdot))     \in L^2(\mathbb{X},\mathbb{R})^M \simeq L^2(\mathbb{X},\mathbb{R}^M).
\end{align*}
Then, we define the approximation $\widehat{K} = \widehat{K}_f^t$ of the Koopman operator~$\mathcal{K} = \mathcal{K}_f^t$ over $\mathbb{V}$  on the compact set $\mathbb{X}$ as a best-fit solution in an $L^2$-sense by means of the regression problem
\begin{equation}\label{eq:fitting}
	\widehat{K} = \argmin_{K \in \mathbb{R}^{M \times M} \hspace*{0.25cm}}\hspace*{-0.2cm} \int_\mathbb{X} \|\Psi(\Flow_f^t(x)) - K\Psi(x)\|^2_2 \,\mathrm{d}x.
\end{equation}
The optimality conditions of \eqref{eq:fitting} read: For all $\Phi \in \mathbb{V}^M$, 
\begin{align*}
	0 = \int_\mathbb{X} \langle \Psi(\Flow_f^t(x)) - \widehat{K} \Psi(x), \Phi(x)\rangle \,\mathrm{d}x = \int_\mathbb{X} \langle \Psi(\Flow_f^t(x)), \Phi(x)\rangle\,\mathrm{d}x - \int_\mathbb{X} \langle \widehat{K} \Psi(x),\Phi(x)\rangle \,\mathrm{d}x
\end{align*}
such that, using $\psi_ie_i\in \mathbb{V}^M$, $i\in [1:M]$, as test functions, where $e_i\in \R^M$ is the $i$-th unit vector, the solution of Problem~\eqref{eq:fitting} can easily be calculated by 
\begin{align}\label{eq:KI}
	\widehat{K} = \Psi_Y \Psi_X^{-1},
\end{align}
where $\Psi_X$ and $\Psi_Y$ are given by $\Psi_X = \int_\mathbb{X} \Psi(x)\Psi(x)^\top\,\mathrm{d}x$ and $\Psi_Y = \int_\mathbb{X}\Psi(\Flow_f^t(x))\Psi(x)^\top \,\mathrm{d}x$, respectively. Here, invertibility of the mass matrix $\Psi_X$ follows from the linear independence of the observables, see~\cite[Lemma C.2]{philipp2024extended}.
The matrix $\widehat{K}$ is a representation of the compression of $\mathcal{K}$, i.e., $\widehat{K} = {P}_\mathbb{V} \mathcal{K} \vert_\mathbb{V}$, where $P_\mathbb{V}$ is the $L^2$-orthogonal projection onto $\mathbb{V}$. 
If the dictionary~$\mathbb{V}$ is Koopman-invariant, then $\widehat{K} = \mathcal{K} \vert_{\mathbb{V}}$ holds. 
If this property is not satisfied, the projection~$P_\mathbb{V}$ introduces a projection error, which was, to the best of our knowledge, first analysed in~\cite{ZhanZuaz23} by means of a dictionary of finite elements, see also~\cite{SchaWort22} for an extension to control systems.
Using a kernel-based dictionary and the rich theory of reproducing kernel Hilbert spaces (RKHS), bounds on the projection error were provided in~\cite{PhilScha23}. 
Recently, also pointwise error bounds in the context of regression in the RKHS norm (instead of the $L^2$-norm) were given in~\cite{KohnPhil24} and also in \cite{yadav2025approximation} using a variant of EDMD with Bernstein polynomials.

Besides the above-mentioned projection error, in data-driven surrogates (corresponding to approximations of the integral by a sum over finitely many data points), an additional estimation error is introduced. We briefly comment on this topic; however, we will keep the idealized situation \eqref{eq:fitting} to illustrate the main contribution of this work. 
\begin{remark}\label{rem:finite_data}
	In fact, the regression problem~\eqref{eq:fitting} is the ideal formulation for the compression in the infinite-data limit.
	An empirical estimator of $\widehat{K}$ given $\new{N}\in \mathbb{N}$ i.i.d.\ data points $x_1,\ldots,x_{\new{N}}$ can be computed via
	$$
	\widehat{K}_{f,\new{N}} = \operatorname{arg} \hspace*{-0.5cm}\min_{K\in \mathbb{R}^{M \times M} \hspace*{0.25cm}}\hspace*{-0.2cm} \sum_{j=1}^{\new{N}} \|\Psi(\Flow_f^t(x_j)) - K \Psi(x_j)\|_2^2.
	$$
	For $\new{N} \geq M$, $\widehat{K}_{f,\new{N}} = (\Psi_X^{\new{N}} {\Psi_Y^{\new{N}}}^\top)(\Psi_X^{\new{N}} {\Psi_X^{\new{N}}}^\top)^{-1}$ is a closed-form solution using the $(M \times \new{N})$-data matrices $\Psi_X^{\new{N}}$ and $\Psi_Y^{\new{N}}$ with entries $\psi_k(x_j)$ and $\psi_k(\Flow_f^t(x_j))$, $(k,j) \in [1:M] \times [1:\new{N}]$, respectively.
	The convergence $\widehat{K}_{f,\new{N}} \rightarrow \widehat{K}$ for $\new{N} \rightarrow \infty$ follows by the law of large numbers~\cite{KordMezi18,Mezi22}.
	For a quantitative convergence analysis and finite-data error bounds we refer to~\cite{NuskPeit23,kostic2024consistent,philipp2024extended} for ergodic sampling, also embedded in the general framework of operator learning, see~\cite{mollenhauer2022kernel} and \cite{colbrook2024rigorous,kostic2024sharp} for spectral approximation properties. 
	For recent results on the estimation error in reproducing kernel Hilbert spaces, we refer to~\cite{PhilScha23}. 
\end{remark}

\section{Reprojections to ensure consistency}\label{sec:problem_formulation}

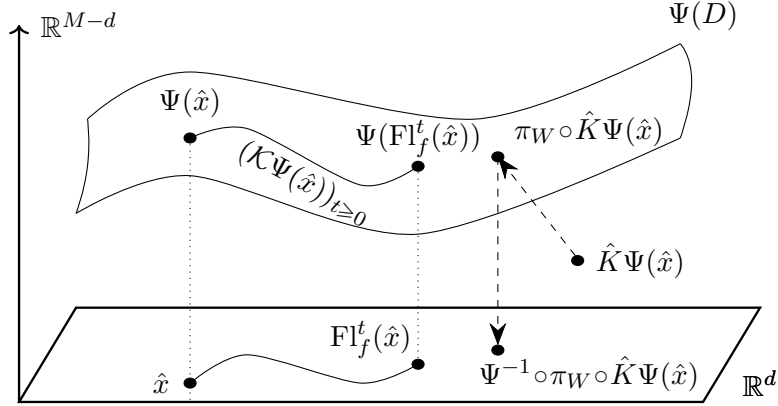
\begin{figure}[htb]
	\centering
	\begin{tikzpicture}[xscale=1.5,yscale=1.25]
    
    \draw[smooth, tension=0.5] plot coordinates{(1,2) (2,2.4) (4,1.78) (6.3,2.8)}{};
    \draw[smooth, tension=0.5] plot coordinates{(6.3,3.8) (4.5,3) (2,3.5) (1.1,3)}{};
    \draw[smooth, tension=1] plot coordinates{(1,2) (1.1,2.5) (1.1,3)}{};
    \draw[smooth, tension=1] plot coordinates{(6.3,2.8) (6.4,3.4) (6.3,3.8)}{};
    \node at (6.5,4.1) {$\Psi(D)$};
    \draw[line width=.08em] (0.5,0)--(1,1)--(7,1)--(6.5,0) -- cycle;
    \node (a) at (7,0.2) {$\mathbb{R}^d$};
    
    \draw[line width=.08em,->] (0.5,0) -- (0.5, 4) node [label=right:$\mathbb{R}^{M-d}$] {};

    \draw[smooth, tension=0.5] plot coordinates{(2,0.2) (2.5,0.5) (3.5,0.2) (4,0.4)} node[above left=-.05cm] {$\Flow_f^t(\hat{x})$};
    \draw[thin, fill=black] (2,0.2) circle (1.5pt) node[label=left:$\hat{x}$]{};
    \draw[thin, fill=black] (4,0.4) circle (1.5pt);

    \draw[smooth, tension=0.5] plot coordinates{(2,2.8) (2.5,2.9) (3.5,2.3) (4,2.5)} node[above=-.01cm] {$\Psi(\Flow_f^t(\hat{x}))$};
    \draw[thin, fill=black] (2,2.8) circle (1.5pt) node[label=above:$\Psi(\hat{x})$]{};
    \draw[thin, fill=black] (4,2.5) circle (1.5pt);
    \draw [dotted] (2,2.8) -- (2,0);
    \draw [dotted] (4,0.4) -- (4,2.5);
    \node (b) at (3,2.32) [rotate=-28] {$(\mathcal{K}\Psi(\hat{x}))_{t\geq 0}$};
    
    \draw[thin, fill = black] (5.4,1.5) circle (1.5pt);
    \node (a) at (5.95,1.5) {$\hat{K}\Psi(\hat{x})$};
    \draw [dashed, -{Stealth[length=3mm]}] (5.38,1.55) -- (4.7,2.6) ;
    \draw[thin, fill=black] (4.7,2.6) circle (1.5pt);
    \node (a) at (7,0.2) {$\mathbb{R}^d$};
    \node (b) at (5.5,2.9) {$\pi_W\!\circ\!\hat{K}\Psi(\hat{x})$};
    \draw [dashed, -{Stealth[length=3mm]}] (4.7,2.6) -- (4.7,0.6);
    \node (c) at (5.5,0.3){$\Psi^{-1}\!\circ\!\pi_W\!\circ\!\hat{K}\Psi(\hat{x})$};
    \draw[thin, fill=black] (4.7,0.55) circle (1.5pt);
\end{tikzpicture}
	\caption{Geometric projection after applying the approximation~$\widehat{K}_f$.}
	\label{fig:manifold}
\end{figure}

\noindent Using the notation of Section~\ref{sec:edmd}, we define the immersed state manifold %
\begin{align}\label{eq:M}
	\mathbb{M} := \im(\Psi) = \{ \Psi(x) \mid x \in \mathbb{X} \} \subset \mathbb{R}^M.
\end{align}
By definition, the set~$\mathbb{M}$ is \textit{invariant} with respect to the Koopman operator~$\mathcal{K}_f^t$, $t \in I = [0,\infty)$, in the sense that %
\begin{equation}\label{eq:Koopman_manifold}
	(\mathcal{K}_f^t \Psi)(x) = \Psi(\Flow_f^t(x)) \in \mathbb{M} \qquad\forall\,x \in \mathbb{X}
\end{equation}
holds, that is, the Koopman operator maps an $\mathbb{M}$-valued function to an $\mathbb{M}$-valued function. 
To fully preserve this property for the compression, see, e.g., \cite{BrunBrun16} for an illustrative example and~\cite{vanGMaho23} for a more detailed discussion, Koopman invariance of the dictionary~$\mathbb{V}$, i.e., if $\mathcal{K} \mathbb{V} \subseteq \mathbb{V}$, is needed since we then have $\widehat{K} = P_\mathbb{V}\mathcal{K}_{\vert \mathbb{V}} =\mathcal{K}_{\vert \mathbb{V}}$. 
Clearly, finding such a Koopman-invariant and still sufficiently expressive dictionary is of high interest as it can be inferred from~\cite{TakeKawa17} and the follow-up works, see, e.g., \cite{HaseCort21} and the references therein. Particularly interesting approaches are, among others, the development of a so-called consistency index in~\cite{HaseCort22} and the use of neural networks in~\cite{GuoKord25}, which is also applicable for parameter-dependent dynamical systems.
However, in general, one cannot expect %
$\widehat{K}\Psi(x) \in \mathbb{M}$ for the approximated Koopman operator~$\widehat{K}$ as depicted in Figure~\ref{fig:manifold}. This becomes even more prominent, if sensors serve as observables as it is typically the case in real-world systems.

\noindent This causes two issues in using $\widehat{K}$ (or data-driven approximations thereof) to model the flow of $f$: 
\begin{enumerate}
	\item The first is that it is unclear how to recover the state values underlying the propagated observables.
	Specifically, if $\widehat{K}\Psi(x) \notin \mathbb{M}$, 
	then by definition there is no value $x \in \mathbb{X} \subset \mathbb{R}^d$ satisfying $\Psi (x) = \widehat{K} \Psi(x)$.
	\item The second issue is that the learning process --~based on the regression problem~\eqref{eq:fitting}~-- only uses measurements of the form $z = \Psi(x)$, i.e., only points contained in the set~$\mathbb{M}$ are taken into account.
	Hence, one cannot expect $\widehat{K} z$ to be meaningful if $z \notin \mathbb{M}$, which may render a repeated application of~$\widehat{K}$ questionable. 
	
\end{enumerate}
Both of these issues can be mitigated by \emph{projecting} the dynamics $z^+ = \widehat{K}z$, $z = \Psi(x)$, \emph{back} to the set~$\mathbb{M}$ after each application of the Koopman operator, see, e.g., \cite{MaurGonc16} and the follow-up work~\cite{MaurGonc19}.
In this paper, we propose a framework for optimal choices of such reprojection maps including an extension to parameter-dependent systems.
In particular, we discuss the following two choices: the widely used \emph{coordinate projection} and our newly proposed \emph{geometric projection}.
Such a projection step is particularly crucial for future applications in Koopman-based (predictive) control using \mbox{EDMDc}~\cite{KordMezi18b} or a bilinear surrogate model~\cite{WillHema2016,BoldGrun23}, where the construction of Koopman-invariant subspaces is a highly-nontrivial issue, see, e.g., \cite{GoswPale21}.

In the conference paper \cite{vanGMaho23}, we formulated first steps towards reprojections for parameter-independent systems. Therein, we considered closest-point projections in weighted norms of the following type.
\begin{definition}\label{def:repr}
	For a positive-definite matrix $W \in \R^{M\times M}$, the \emph{weighted reprojection} $\pi_W : \R^{M} \to \mathbb{M}$ satisfies 
	\begin{align}\label{eq:reprojection_function}
		\pi_W(z) \in \operatorname{arg}\min_{q \in \mathbb{M}} \vert q - z \vert_W^2
	\end{align}
	with the weighted norm $|z|_W^2 := z^\top W z$. %
	That is, for a point $z \in \R^{M}$, $\pi_W(z)$ provides the closest points in $\mathbb{M}$ with respect to the norm induced by $W$. %
\end{definition}
Hence, after pursuing one or several steps along the data-driven model, we may reproject to the manifold $\mathbb{M}$ by means of such a reprojection. In this work, we provide a thorough analysis of this idea addressing the following points:
\begin{enumerate}
	\item Tailored application to parameter- or control-affine systems extending the results of~\cite{vanGMaho23} for autonomous systems.
	\item Sophisticated data-driven choice of the weight matrix~$W$ by means of covariance matrices, allowing for a maximum likelihood interpretation of the projection~$\pi_W$.
	\item Efficient evaluation of~\eqref{eq:reprojection_function} using Riemannian Newton methods leveraging the structure of the manifold~$\mathbb{M}$.
\end{enumerate}

\new{We stress that any reprojection scheme entails a trade-off, which, to the best of our knowledge, is rarely made explicit in the literature, where reprojections often only appear implicitly in the implementation. A main advantage of Koopman-based surrogates is their \emph{global linearity}, i.e., a single matrix $\widehat{K}$ propagates the lifted state over arbitrarily long horizons, rendering the well-established machinery of linear systems theory, spectral analysis and linear predictive control directly applicable. Interleaving this linear propagation with the nonlinear map~$\pi_W$ clearly destroys this property, as the resulting predictor $\pi_W \circ \widehat{K}$ is nonlinear on $\mathbb{M}$ and its long-term behaviour may no longer be inferred from the spectrum of $\widehat{K}$ alone.
	However, global linearity is only an advantage as long as the linear model is a faithful surrogate of the dynamics. As the numerical examples in Subsection~\ref{subsec:numres} show, this is usually not the case beyond a few prediction steps, since the iterates without reprojection leave $\mathbb{M}$, become inconsistent with the lifting and, in the examples considered here, either converge to spurious equilibria or diverge. Hence, we accept this loss in order to gain consistency with the lifting and, thus, reliable long-term predictions.
	Two aspects mitigate this loss. First, the model remains exactly linear between two reprojections such that the linear machinery is available on the respective prediction windows. The adaptive scheme of Subsection~\ref{subsec:multistep} renders these windows as long as the desired accuracy permits and, thus, allows to tune this trade-off. Second, if the dictionary spans a Koopman-invariant subspace, the reprojection reduces to the identity on $\mathbb{M}$ and nothing is lost. Hereby, the reprojection may be interpreted as a correction whose magnitude quantifies the deviation from Koopman invariance.}
This reprojection is optimal for recovering the state of the system in the sense that, if we take the model %
$\hat K \Psi(x)$ to be correct, then the most likely admissible state of the system is given by the reprojection $\pi_W(\hat K \Psi(x))$ upon a suitable choice of the weighting matrix $W$.
\begin{proposition}\label{prop:maxlikely}
	Let $x \in \mathbb{X}$ and $z = \Psi(x) + \mu \new{\in} \R^M$, where $\mu \sim N(0, \Sigma)$ for a positive definite covariance $\Sigma$.
	Then the maximum likelihood estimate $\hat{z} \in \mathbb{M}$ of $\Psi(x)$ given a measurement~$z$ satisfies
	\begin{align*}
		\hat{z} \in \pi_{\Sigma^{-1}}(z),
	\end{align*}
	i.e., it is a weighted reprojection~\eqref{eq:reprojection_function} with weighting matrix $\new{W = \Sigma^{-1}}$.
\end{proposition}
\begin{proof}
	Up to a constant offset, the negative log-likelihood of $\hat{z} = \Psi(x)$ given $z$ is
	\begin{align*}
		l(\hat{z} = \Psi(x) \;|\; z) = \frac{1}{2} \vert z - \hat{z} \vert^2_{\Sigma^{-1}}.
	\end{align*}
	An optimal solution for $\hat{z}$ therefore satisfies $\hat{z} = \pi_{\Sigma^{-1}}(z)$.
\end{proof}
\new{We stress the distinction between the two matrices appearing in Proposition~\ref{prop:maxlikely}: $\Sigma$ is the covariance of the perturbation~$\mu$, whereas $W$ is the weighting matrix inducing the metric in~\eqref{eq:reprojection_function}. The maximum-likelihood reprojection corresponds to the choice $W = \Sigma^{-1}$, i.e., directions of large uncertainty are penalized weakly and directions of small uncertainty are penalized strongly. Throughout this work, $W$ denotes a metric, i.e., an inverse covariance, and never a covariance itself.}
A widely-used and extremely simple form is the coordinate projection that was introduced, to the best of our knowledge, in~\cite{MaurGonc16}.
\begin{remark}\label{rem:coord}
	Assume that the first $d$ entries of the dictionary are the coordinate functions, i.e., $\Psi_{i}(x) = x_i$ for $i\in [1:d]$. Then the weight function %
	\begin{align*}
		W_c = \begin{pmatrix}
			I_d & 0\\
			0&0
		\end{pmatrix}\in \R^{M\times M}
	\end{align*}
	induces a coordinate reprojection, that is,
	\begin{align*}
		\Psi(z_{[1:d]}) =\pi_{W_c}(z) \in \operatorname{arg}\min_{q \in \mathbb{M}} \vert q - z \vert_{W_c}^2 ,
	\end{align*}
	where $z_{[1:d]} \in \R^d$ corresponds to the first $d$ entries of $z\in \R^M$. For more details, i.e., a proof that $W_c$ indeed induces a Riemannian metric on $\mathbb{M}$, we refer the reader to \cite[Proposition 5.2]{vanGMaho23}.
\end{remark}

\section{Reprojections for Parameter-Affine Systems}
\label{sec:parameter-affine}

\noindent We consider parameter-affine systems of the form
\begin{align}\label{eq:system_dfn}
	\dot{x} = f(x) + \sum_{i=1}^m p_i g_i(x),
\end{align}
where $f, g_1,...,g_m : \mathbb{X} \to \R^d$ are locally Lipschitz continuous vector fields and $p \in \mathbb{R}^m$, $m\in \mathbb{N}$, represents a constant parameter. An extension to parameter-varying problems (e.g.~control systems) is straightforward such that in this work, we stick to the constant setting for ease of notation. In the following, for $p\in \R^m$ and $t\geq 0$, we denote by $\mathcal{K}_{f + \sum_{i=1}^m p_i g_i}^t$ the Koopman operator induced by $\eqref{eq:system_dfn}$, that is, for $x\in \mathbb{X}$ and $\varphi:\mathbb{X}\to \R$,
\begin{align*}
	(\mathcal{K}_{f + \sum_{i=1}^m p_i g_i}^t\varphi)(x) = \varphi(\Flow_{f + \sum_{i=1}^m p_i g_i}^t(x)),
\end{align*}
again tacitly assuming that $\mathbb{X}$ is invariant under $ \Flow^t_{f + \sum_{i=1}^m p_i g_i}$ for simplicity. 
Since one typically restricts parameter ranges to relatively small intervals, see, e.g., the examples below in Subsections~\ref{subsec:duffing} and \ref{subsec:lorenz}, this is a reasonably mild assumption for dynamical systems.

\subsection{The Koopman operator for parameter-affine systems}
\label{subsec:parameter-affine}

While the parameter-affine structure is directly inherited by the Koopman generator~\cite{PeitOtto20,Sura16,WillHema2016}, see also \cite{guo2025modularized} where this was leveraged for efficient learning on graphs, the Koopman operator is only approximately parameter-affine, as may be seen with a Taylor-series argument: For an observable $\varphi \in C^2(\mathbb{X}, \R)$ and any $x \in \mathbb{X}$, one has
\begin{align*}
	(\mathcal{K}_{f + \sum_{i=1}^m p_i g_i}^t \varphi) (x) & = \varphi (\Flow_{f + \sum_{i=1}^{m} p_i g_i}^t (x)) \\
	& = \varphi(x) + t \underbrace{\tD \varphi (x)(f + \sum_{i=1}^m p_i g_i)(x)}_{= \tD \varphi (x)f(x) + \sum_{i=1}^m p_i \tD \varphi (x) g_i(x)} +\ \mathcal{O}(t^2) \\
	& = \Big[ \varphi(x) + t \tD \varphi (x)f(x) \Big] + \sum_{i=1}^m p_i \Big[ t \tD \varphi (x) g_i(x) \pm \varphi(x) \Big] + \mathcal{O}(t^2) \\
	&= \varphi(\Flow_f^t(x)) + \sum_{i=1}^m p_i  \Big[ \varphi(\Flow_{g_i}^t (x)) -\varphi(x) \Big] + \mathcal{O}(t^2)\\
	&= (\mathcal{K}^t_f\varphi)(x) + \sum_{i=1}^m p_i \left( \mathcal{K}^t_{g_i} \varphi - \varphi\right)(x) + \mathcal{O}(t^2)
\end{align*}
Motivated by the above calculation and denoting by $\otimes$ the usual Kronecker product we obtain the (infinite-dimensional) approximation for $\Psi$ as defined in \eqref{eq:Psi} 
\begin{align}\label{eq:firstapprox}
	(\mathcal{K}_{f + \sum_{i=1}^m p_i g_i}^t \Psi)(x)
	&= \widebar{\mathcal{K}}^t (\bar{p} \otimes\Psi)(x) + \mathcal{O}(t^2),
\end{align}
with the augmented parameter vector $\bar{p}:=\left(\begin{matrix}1\\p\end{matrix}\right)\in \R^{m+1}$  %
and the column block operator
\begin{align}\label{eq:bilinear1}
	\widebar{\mathcal{K}}^t := \begin{pmatrix}
		\mathcal{K}_f^t &
		\mathcal{K}_{g_1}^t - \mathcal{I} & \cdots &
		\mathcal{K}_{g_m}^t - \mathcal{I}
	\end{pmatrix},
\end{align}
where $\mathcal{K}^t_h$,  $h\in \{f,g_1,\ldots,g_m\}$, is the Koopman operator corresponding to the flow $\Flow^t_h$ defined by~\eqref{eq:Koopman}.

This bilinear representation motivates the following (finite-dimensional) approximation
\begin{align}\label{eq:parameter_koopman}
	\widehat{K} &:=
	\operatorname{argmin}_{K \in \mathbb{R}^{M \times M(m+1)} \hspace*{0.4cm}}\hspace*{-0.2cm} \int_{\mathbb{L} \otimes \mathbb{X}} \|\Psi(\Flow_{f + \sum_{i=1}^m p_i g_i}^t(x)) - {K}(\bar{p} \otimes \Psi(x))\|^2_2 \,\mathrm{d}(p \otimes x),
\end{align}
where $\mathbb{L}\subset \R^m$ is a compact set (with potentially only finitely-many elements) containing all possible parameters~$p$.
Solving the problem~\eqref{eq:parameter_koopman} yields an $M \times M(m+1)$ matrix that can be partitioned into 
\begin{align}\label{eq:blockcols}
	\widehat{K} = \begin{pmatrix}
		\widehat{K}_0 &
		\widehat{K}_1 & \cdots &
		\widehat{K}_{m}
	\end{pmatrix}
\end{align}
with $\widehat{K}_i \in \R^{M \times M}$ for $i \in [0:m]$, from which the individual approximations of the block-columns of $\widebar{\mathcal{K}}^t$ in~\eqref{eq:bilinear1} can be defined. This leads to the approximation of the parametric Koopman operator
\begin{align}\label{eq:parameter_kooopman_model}
	\left( \mathcal{K}^t_{f+\sum_{i=1}^m p_i g_i} \Psi \right)(x) = \Psi(\Flow^t_{f + \sum_{i=1}^m p_i g_i}(x))
	\approx \widehat{K}_0 \Psi(x) + \sum_{i=1}^m p_i \widehat{K}_i \Psi(x).
\end{align}
We stress that the structure of the right-hand side of \eqref{eq:bilinear1} mainly serves as a motivation to formulate the regression problem \eqref{eq:parameter_koopman} that may be leveraged in an error analysis. %

As in the parameter-independent case \eqref{eq:Koopman_manifold}, the parametric Koopman operator on the left-hand side of \eqref{eq:parameter_kooopman_model} preserves the set $\mathbb{M} = \Psi(\mathbb{X})$.
However, as in the parameter-independent case, this does, in general, not hold for the approximation on the right-hand side of~\eqref{eq:parameter_kooopman_model}. %
Thus, the state of the approximating system must be reprojected to~$\mathbb{M}$ to recover a state of the original system in $\mathbb{X}$. Before presenting a suitable reprojection in Section~\ref{subsec:reprojections}, we briefly comment on data-driven surrogates of \eqref{eq:parameter_koopman} via EDMD and provide a relation to existing parameter-affine approximations of the Koopman operator.
\begin{remark}\hfill
	\begin{itemize}
		\item[(i)] 
		In the spirit of EDMD and analogously to Remark~\ref{rem:finite_data} we may perform a Monte-Carlo-like approximation of \eqref{eq:parameter_koopman} by replacing the integral by a sum over data points, e.g., chosen i.i.d.\ with respect to the Lebesgue measure on $\mathbb{L}\otimes \mathbb{X}$. 
		\item[(ii)] 	Following \cite[Section 5]{PhilScha23}, again using a Taylor-series argument, we have 
		\begin{align*}
			\left( \mathcal{K}^t_{g_i} \varphi - \varphi \right)(x) %
			& = t\nabla \varphi \cdot (f + g_i - f)(x) %
			+ \mathcal{O}(t^2)\\
			& = \left(\mathcal{K}^t_{f+g_i} \varphi - \varphi\right)(x) - \left(\mathcal{K}^t_f \varphi - \varphi\right)(x) + \mathcal{O}(t^2)\\
			& = \left( \mathcal{K}^t_{f+g_i} \varphi \right)(x) - \left( \mathcal{K}^t_f\varphi \right)(x) + \mathcal{O}(t^2),
		\end{align*}
		which shows the equation
		\begin{align*}
			\mathcal{K}^t_f\varphi (\hat x) + \sum_{i=1}^m p_i \left( \mathcal{K}^t_{g_i} \varphi - \varphi\right) (\hat x) = \mathcal{K}^t_f\varphi(\hat x) + \sum_{i=1}^m p_i \left( \mathcal{K}^t_{f + g_i} \varphi - \mathcal{K}^t_f\varphi\right) (\hat x) + \mathcal{O}(t^2).
		\end{align*}
		The presented line of reasoning allows us to replace the definition ~\eqref{eq:bilinear1} by %
		\begin{align}\label{eq:bilinear2}
			\bar{\mathcal{K}}^t := \begin{pmatrix}
				\mathcal{K}_f^t & \mathcal{K}_{g_1+f}^t - \mathcal{K}_f^t & \cdots & \mathcal{K}_{g_m+f}^t - \mathcal{K}_f^t
			\end{pmatrix},
		\end{align}
		which nicely links the presented parameter-affine approximation of the Koopman operator to existing ones~\cite{PeitOtto20,NuskPeit23}. 
	\end{itemize}
\end{remark}

\subsection{Maximum likelihood reprojections via covariance estimators}\label{subsec:reprojections}

\noindent Consider the parameter-affine system \eqref{eq:system_dfn} and let $\mathbb{V}$ be a $M$-dimensional dictionary.
Let $\Psi \in L^2(\mathbb{X}, \R)^M$ denote the vector of observables spanning~$\mathbb{V}$ and $\widehat{K}$ be the approximation of the Koopman operator defined in~\eqref{eq:parameter_koopman}.
Then, as discussed in Section~\ref{sec:problem_formulation}, it is generally not guaranteed that the set {$\mathbb{M} = \Psi(\mathbb{X})$} is preserved under the \karl{approximated Koopman operator~$\widehat{K}$}, {compare~\eqref{eq:Koopman_manifold}}. %
To measure the approximation error, %
the \emph{residual operator} for a vector field $f$ is defined on $\mathbb{V}$ by
\begin{align*}
	(\mathcal{R}_{f}^t \Psi)(x) := (\mathcal{K}_{f}^t \Psi)(x) - {\widehat{K}_{f}} \Psi(x),
\end{align*}
where by $\widehat{K}_f$ we denote the solution of \eqref{eq:fitting}.

In the following, we model the residual as a zero-mean random variable. To this end, we impose a mild condition on the dictionary, which we tacitly assume from now on.
\begin{assumption}[\new{Constant observable}]\label{ass:constant}
	The constant function is contained in the dictionary, i.e., $\one \in \mathbb{V}$, where $\one(x) = 1$ for all $x \in \mathbb{X}$. Without loss of generality, we set $\psi_1 = \one$.
\end{assumption}
Assumption~\ref{ass:constant} holds for every polynomial or monomial dictionary containing the monomial of order zero and, in particular, for all dictionaries used in Subsection~\ref{subsec:numres}. It is imposed for the following reason. Written in matrix form, the optimality conditions of the regression problem~\eqref{eq:fitting} stated in Section~\ref{sec:edmd} read
\begin{align}\label{eq:normalequations}
	\int_{\mathbb{X}} (\mathcal{R}^t_f \Psi)(x)\, \Psi(x)^\top \,\mathrm{d}x = 0 \in \R^{M \times M},
\end{align}
i.e., the residual is orthogonal to the span of the dictionary, but its mean does not necessarily vanish. However, under Assumption~\ref{ass:constant}, the first column of~\eqref{eq:normalequations} equals $\int_{\mathbb{X}} (\mathcal{R}^t_f \Psi)(x)\,\mathrm{d}x$ such that the residual is centered. Note that, if Assumption~\ref{ass:constant} is dropped, all subsequent statements remain valid upon replacing $\mathcal{R}^t_f \Psi$ by the centered residual $\mathcal{R}^t_f \Psi - \mathbb{E}(\mathcal{R}^t_f\Psi)$ and adding the respective constant offset to the prediction model.

Then, by virtue of the definition of $\widehat{K}$ through the regression problem \eqref{eq:fitting}, Assumption~\ref{ass:constant} and~\eqref{eq:normalequations}, the first and second moments of $\mathcal{R}^t_f \Psi$ satisfy
\begin{align*}
	\mathbb{E}(\mathcal{R}^t_f \Psi) & = \int_{\mathbb{X}} \mathcal{R}^t_f \Psi(x) \mathrm{d} x = \int_{\mathbb{X}} \Psi(\Flow_f^t(x)) - \widehat{K}_{\new{f}} \Psi(x)\, \mathrm{d} x = 0, \\
	\mathrm{Cov}(\mathcal{R}^t_f \Psi)
	&= \mathbb{E} \left((\mathcal{R}^t_f \Psi(x) - \mathbb{E}(\mathcal{R}^t_f \Psi))(\mathcal{R}^t_f \Psi(x) - \mathbb{E}(\mathcal{R}^t_f \Psi))^\top\right) \\
	&= \int_{\mathbb{X}} \mathcal{R}^t_f \Psi(x)(\mathcal{R}^t_f \Psi(x))^\top \mathrm{d} x \\
	&= \int_{\mathbb{X}} (\Psi(\Flow_f^t(x)) - \widehat{K}_{\new{f}} \Psi(x))(\Psi(\Flow_f^t(x)) - \widehat{K}_{\new{f}} \Psi(x))^\top \mathrm{d} x.
\end{align*}
Here, we propose a model for the approximated Koopman dynamics as a linear system corrupted by a zero-mean Gaussian noise process, that is, %
\begin{align*}
	\Psi(\Flow_f^t(x)) &= \widehat{K}_f \Psi(x) + \mathcal{R}_f^t \Psi(x)
	\approx \widehat{K}_{\new{f}} \Psi(x) + \mu, &
	\mu &\sim N(0, \mathrm{Cov}(\mathcal{R}^t_f \Psi)).
\end{align*}
While this is also an approximation, it captures the effects of the finite-dictionary Koopman fitting up to second order.
Moreover, it provides a stochastic approximation of the Koopman dynamics %
including a model of the error introduced by projection.
In Proposition \ref{prop:cov}, we demonstrate how the covariance used to model the error can be computed from the same data used to construct $\widehat{K}$ in the first place.

For parameter-dependent systems, we define the residual operator analogously, that is,
\begin{align}\label{eq:residual_operator}
	\mathcal{R}^t_{f + \sum_{i=1}^m p_i g_i} \Psi (x) :=
	(\mathcal{K}^t_{f + \sum_{i=1}^m p_i g_i} \Psi)(x) -  \left(\widehat{K}_0 + \sum_{i=1}^m p_i \widehat{K}_i\right) \Psi(x),
\end{align}
for all initial states $x \in \mathbb{X}$ and parameters $p \in \mathbb{L}$. Here, we recall that $\widehat{K}_i\in \mathbb{R}^{M\times M}$, $i\in [0:m]$, are defined as the block columns of $\widehat{K}$, see \eqref{eq:blockcols}, solving the augmented regression problem~\eqref{eq:parameter_koopman}.
\new{In contrast to the autonomous case, the regression problem~\eqref{eq:parameter_koopman} is posed on the product $\mathbb{L}\otimes\mathbb{X}$ such that its optimality conditions only enforce a vanishing mean with respect to the \emph{joint} measure. Since the reprojection is applied for a \emph{fixed} parameter~$p$, we make the respective statement precise in the following.}
\begin{lemma}[\new{Centered residual at fixed parameter}]\label{lem:centered}
	\new{Let Assumption~\ref{ass:constant} hold and assume that $\mathbb{L}$ is not contained in an affine hyperplane of $\R^m$ up to a set of measure zero, i.e., that the moment matrix $G := \int_{\mathbb{L}} \bar{p}\bar{p}^\top \,\mathrm{d}p \in \R^{(m+1)\times(m+1)}$ is invertible. Then}
	\begin{align*}
		\new{\mathbb{E}\big(\mathcal{R}^t_{f + \sum_{i=1}^m p_i g_i} \Psi\big) = \mathcal{O}(t^2)
			\qquad \text{uniformly in } p \in \mathbb{L}.}
	\end{align*}
\end{lemma}
\begin{proof}
	\new{Abbreviate $R(\cdot,p) := \mathcal{R}^t_{f + \sum_{i=1}^m p_i g_i} \Psi$. Inserting the Taylor expansion of Subsection~\ref{subsec:parameter-affine} into~\eqref{eq:residual_operator} and collecting the terms of order zero and one in $p$ yields}
	\begin{align}\label{eq:affineresidual}
		\new{R(x,p) = r_0(x) + \sum_{i=1}^m p_i\, r_i(x) + \mathcal{O}(t^2),}
	\end{align}
	\new{where $r_0 := \mathcal{K}^t_f \Psi - \widehat{K}_0 \Psi$ and $r_i := (\mathcal{K}^t_{g_i} \Psi - \Psi) - \widehat{K}_i \Psi$ for $i \in [1:m]$. Note that these functions are defined exclusively by means of the block columns~$\widehat{K}_i$ of the augmented regression problem~\eqref{eq:parameter_koopman}. Hence, $R(\cdot,p)$ is affine in~$p$ up to terms of order $\mathcal{O}(t^2)$. The optimality conditions of~\eqref{eq:parameter_koopman} read $\int_{\mathbb{L}} \int_{\mathbb{X}} R(x,p) (\bar p \otimes \Psi(x))^\top \,\mathrm{d}x\,\mathrm{d}p = 0$. Thus, under Assumption~\ref{ass:constant}, evaluating the block columns corresponding to $\psi_1 = \one$ yields}
	\begin{align*}
		\new{\int_{\mathbb{L}} \bar{p}_i\, q(p) \,\mathrm{d}p = 0 \quad \text{for all } i \in [0:m],
			\qquad \text{where } q(p) := \int_{\mathbb{X}} R(x,p) \,\mathrm{d}x \in \R^M.}
	\end{align*}
	\new{By~\eqref{eq:affineresidual}, we have $q(p) = \bar{R}\,\bar{p} + \mathcal{O}(t^2)$ with $\bar{R} := \left(\int_\mathbb{X} r_0 \,\mathrm{d}x \;\cdots\; \int_\mathbb{X} r_m \,\mathrm{d}x\right) \in \R^{M \times (m+1)}$ such that the above conditions read $\bar{R}\,G = \mathcal{O}(t^2)$. Hence, invertibility of~$G$ yields $\bar{R} = \mathcal{O}(t^2)$ and, thus, $q(p) = \mathcal{O}(t^2)$ for all $p \in \mathbb{L}$, where uniformity follows from compactness of~$\mathbb{L}$.}
\end{proof}
\new{We stress that, in general, the mean does \emph{not} vanish exactly for a fixed parameter~$p$, but only up to second-order terms, which are neglected in the parameter-affine model~\eqref{eq:parameter_kooopman_model}. This explains the terms of order $\mathcal{O}(t^2)$ occurring in the proof of Proposition~\ref{prop:cov} below. Note that, if $\mathbb{L}$ consists of finitely many parameter samples, invertibility of~$G$ amounts to at least $m+1$ affinely independent samples.}
For parameter-dependent systems~\eqref{eq:system_dfn}, the covariance may depend on the parameter $p$ %
in a non-trivial fashion. However, in the following, we will show that it can be approximated as a quadratic function of the parameters for short integration times. To this end, define the matrix-valued bilinear form $Q : \R^{m+1} \times \R^{m+1} \to \R^{M \times M}$ as the symmetric solution to the regression problem
\begin{align}\label{eq:solQ}
	\min_Q \int_{\mathbb{L}} \left\| 
	\operatorname{Cov}(\mathcal{R}^t_{f + \sum_{i=1}^m p_i g_i} \Psi) - Q(\bar{p}, \bar{p})
	\right\|_F^2 \mathrm{d}p,
\end{align}
where $\bar{p} = \left(\begin{smallmatrix}1\\p\end{smallmatrix}\right)\in \R^{m+1}$ again denotes the augmented %
parameter vector.
The solution $Q$ may be determined analytically, that is, denoting by $\operatorname{vec}$ the column-wise vectorization of a matrix, $Q$ solves~\eqref{eq:solQ} if
\begin{align}\label{eq:covariance_map_Q}
	\operatorname{vec} (Q(\bar{p}, \bar{p})) &= Y X^\dag (\bar{p} \otimes \bar{p}), \\
	\mathrm{where\ }Y &:= \int_{\mathbb{L}} \int_{\mathbb{X}} (\mathcal{R}^t_{f + \sum_{i=1}^m p_i g_i} \Psi(x) \bar{p}^\top ) \otimes (\mathcal{R}^t_{f + \sum_{i=1}^m p_i g_i} \Psi(x) \bar{p}^\top ) 
	\; \mathrm{d}x \; \mathrm{d}p, \notag\\
	X &:= \int_{\mathbb{L}} \bar{p} \bar{p}^\top \otimes \bar{p} \bar{p}^\top \mathrm{d}p. \notag
\end{align}
Hereby we note that $X$ may be computed analytically, where $Y$ depends on the Koopman operator and thus has to be approximated by means of data in applications, see Remark~\ref{rem:edmdcov}.

The following result shows that the bilinear form $Q$ yields a second-order approximation of the true covariance.

\begin{proposition}\label{prop:cov}
	For any $p \in \mathbb{L}$, the covariance of the residual operator \eqref{eq:residual_operator} satisfies
	\begin{align}\label{eq:koopman_system_noise}
		\operatorname{Cov}\left(
		\mathcal{R}^t_{f + \sum_{i=1}^m p_i g_i} \Psi
		\right)
		&= Q(\bar{p},\bar{p}) + \mathcal{O}(t^2),
	\end{align}
	where $Q$ is defined as in \eqref{eq:covariance_map_Q}.
\end{proposition}

\begin{proof}
	For the expected value, we have that %
	\begin{align*}
		\mathbb{E}(\mathcal{R}^t_{f + \sum_{i=1}^m p_i g_i} \Psi)
		&= \mathbb{E}(\mathcal{R}^t_{f}\Psi + \sum_{i=1}^m p_i \mathcal{R}^t_{g_i} \Psi) + \mathcal{O}(t^2) \\
		&= \mathbb{E}(\mathcal{R}^t_{f}\Psi) + \sum_{i=1}^m p_i \mathbb{E}(\mathcal{R}^t_{g_i} \Psi)  + \mathcal{O}(t^2) \\
		&= \mathcal{O}(t^2)
	\end{align*}
	as $\mathbb{E}(\mathcal{R}^t_{f}\Psi) = \mathbb{E}(\mathcal{R}^t_{g_1}\new{\Psi}) = \ldots = \mathbb{E}(\mathcal{R}^t_{g_m}\new{\Psi}) = \mathcal{O}(t^2)$\new{, cf.\ Lemma~\ref{lem:centered}}.
	For convenience of notation, %
	let $g_0 := f$, so that
	\begin{align*}
		f + \sum_{i=1}^m p_i g_i = \sum_{i=0}^m \bar{p}_i g_i.
	\end{align*}
	Then, using the parameter-affine structure, the covariance satisfies %
	\begin{align}
		&\mathrm{Cov}(\mathcal{R}^t_{f + \sum_{i=1}^m p_i g_i} \Psi)\nonumber
		\\&= \mathrm{Cov}(\mathcal{R}^t_{\sum_{i=0}^m \bar{p}_i g_i} \Psi) \nonumber\\
		&= \mathrm{Cov}\left(\sum_{i=0}^m \bar{p}_i \mathcal{R}^t_{g_i} \Psi\right) + \mathcal{O}(t^2)\nonumber \\
		&= \mathbb{E} \left( \left(\sum_{i=0}^m \bar{p}_i \mathcal{R}^t_{g_i} \Psi - \underbrace{\mathbb{E}\left(\sum_{i=0}^m \bar{p}_i \mathcal{R}^t_{g_i} \Psi\right)}_{=\mathcal{O}(t^2)}\right) \left(\sum_{i=0}^m \bar{p}_i \mathcal{R}^t_{g_i} \Psi - \underbrace{\mathbb{E}\left(\sum_{i=0}^m \bar{p}_i \mathcal{R}^t_{g_i} \Psi\right)}_{=\mathcal{O}(t^2)}\right)^\top\right) + \mathcal{O}(t^2)\nonumber\\
		&= \sum_{i,j=0}^m \bar{p}_i \bar{p}_j \mathbb{E} \left( (\mathcal{R}^t_{g_i} \Psi)( \mathcal{R}^t_{g_j} \Psi)^\top \right)  + \mathcal{O}(t^2).\label{eq:almostdone}
	\end{align}
	We now show that $Q$ defined by \eqref{eq:covariance_map_Q} provides an approximation of the first term in \eqref{eq:almostdone} up to second order in $t$. %
	To see this, we compute
	\begin{align*}
		Y &= \int_{\mathbb{L}} \int_{\mathbb{X}}
		(\mathcal{R}^t_{f + \sum_{i=1}^m p_i g_i} \Psi(x) \bar{p}^\top )
		\otimes (\mathcal{R}^t_{f + \sum_{i=1}^m p_i g_i} \Psi(x) \bar{p}^\top ) 
		\; \mathrm{d}x \; \mathrm{d}p \\
		&= \int_{\mathbb{L}} \int_{\mathbb{X}}
		(\mathcal{R}^t_{\sum_{i=0}^m \bar{p}_i g_i} \Psi(x) \bar{p}^\top )
		\otimes (\mathcal{R}^t_{\sum_{i=0}^m \bar{p}_i g_i} \Psi(x) \bar{p}^\top )
		\; \mathrm{d}x \; \mathrm{d}p \\
		&= \int_{\mathbb{L}} \int_{\mathbb{X}}
		\left(\left(\begin{matrix}
			\mathcal{R}^t_{g_0}\Psi(x) & \cdots & \mathcal{R}^t_{g_m}\Psi(x)
		\end{matrix}\right) \bar{p} \bar{p}^\top\right)
		\otimes \left(\left(\begin{matrix}
			\mathcal{R}^t_{g_0}\Psi(x) & \cdots & \mathcal{R}^t_{g_m}\Psi(x)
		\end{matrix}\right) \bar{p} \bar{p}^\top\right)
		\; \mathrm{d}x \; \mathrm{d}p
		+ \mathcal{O}(t^2) \\
		&= \int_{\mathbb{L}} \int_{\mathbb{X}}
		\left(\left(\begin{matrix}
			\mathcal{R}^t_{g_0}\Psi(x) & \cdots & \mathcal{R}^t_{g_m}\Psi(x)
		\end{matrix}\right) 
		\otimes \left(\begin{matrix}
			\mathcal{R}^t_{g_0}\Psi(x) & \cdots & \mathcal{R}^t_{g_m}\Psi(x)
		\end{matrix}\right) \right)
		(\bar{p} \bar{p}^\top \otimes
		\bar{p} \bar{p}^\top)
		\; \mathrm{d}x \; \mathrm{d}p
		+ \mathcal{O}(t^2) \\
		&= \int_{\mathbb{X}}
		\left(\left(\begin{matrix}
			\mathcal{R}^t_{g_0}\Psi(x) & \cdots & \mathcal{R}^t_{g_m}\Psi(x)
		\end{matrix}\right) 
		\otimes \left(\begin{matrix}
			\mathcal{R}^t_{g_0}\Psi(x) & \cdots & \mathcal{R}^t_{g_m}\Psi(x)
		\end{matrix}\right) \right)
		\; \mathrm{d}x 
		\int_{\mathbb{L}} 
		(\bar{p} \bar{p}^\top \otimes \bar{p} \bar{p}^\top)
		\; \mathrm{d}p
		+ \mathcal{O}(t^2).
	\end{align*}
	Substituting this into the definition of $Q$, i.e., \eqref{eq:covariance_map_Q}, we get %
	\begin{align*}
		\operatorname{vec} (Q(\bar{p}, \bar{p})) &= Y X^\dag (\bar{p} \otimes \bar{p}) \\
		&= \int_{\mathbb{X}}
		\left(\left(\begin{matrix}
			\mathcal{R}^t_{g_0}\Psi(x) & \cdots & \mathcal{R}^t_{g_m}\Psi(x)
		\end{matrix}\right) 
		\otimes \left(\begin{matrix}
			\mathcal{R}^t_{g_0}\Psi(x) & \cdots & \mathcal{R}^t_{g_m}\Psi(x)
		\end{matrix}\right) \right)
		\; \mathrm{d}x (\bar{p} \otimes \bar{p}) + \mathcal{O}(t^2) \\
		&= \int_{\mathbb{X}}
		\left(\left(\begin{matrix}
			\mathcal{R}^t_{g_0}\Psi(x) & \cdots & \mathcal{R}^t_{g_m}\Psi(x)
		\end{matrix}\right) \bar{p} \right)
		\otimes \left(\left(\begin{matrix}
			\mathcal{R}^t_{g_0}\Psi(x) & \cdots & \mathcal{R}^t_{g_m}\Psi(x)
		\end{matrix}\right) \bar{p}\right)
		\; \mathrm{d}x + \mathcal{O}(t^2) \\
		&=  \int_{\mathbb{X}}{
			\left( \mathcal{R}^t_{g_0}\Psi(x) p_0 + \cdots + \mathcal{R}^t_{g_m}\Psi(x) p_m \right)}
		\otimes { \left( \mathcal{R}^t_{g_0}\Psi(x) p_0 + \cdots + \mathcal{R}^t_{g_m}\Psi(x) p_m \right)}
		\; \mathrm{d}x + \mathcal{O}(t^2)
	\end{align*}
	such that
	\begin{align*}
		Q(\bar{p}, \bar{p})
		&= \int_{\mathbb{X}} 
		{ \left( \mathcal{R}^t_{g_0}\Psi(x) p_0 + \cdots + \mathcal{R}^t_{g_m}\Psi(x) p_m \right) \left( \mathcal{R}^t_{g_0}\Psi(x) p_0 + \cdots + \mathcal{R}^t_{g_m}\Psi(x) p_m \right)^\top}
		\; \mathrm{d}x + \mathcal{O}(t^2) \\
		&= \sum_{i,j=0}^m \bar{p}_i \bar{p}_j \mathbb{E} \left( (\mathcal{R}^t_{g_i} \Psi)( \mathcal{R}^t_{g_j} \Psi)^\top \right)  + \mathcal{O}(t^2).
	\end{align*}
	This completes the proof.
\end{proof}

Combining Proposition~\ref{prop:maxlikely} and Proposition~\ref{prop:cov}, we directly obtain the following result for the parameter-dependent prediction model.
\begin{corollary}\label{cor:mlparam}
	Let $x \in \mathbb{X}$ and $z = \Psi(x) \in \mathbb{M} \subset \R^M$. Let
	\begin{align*}
		z^+ := \left(\widehat{K}_0 + \sum_{i=1}^m p_i \widehat{K}_i\right) z
		= \Psi(\Flow_{f + \sum_{i=1}^m p_i g_i}^t(x))
		+ \mu, &&
		\mu \sim N\left(0, \sum_{i,j=0}^m Q_{ij} p_i p_j\right).
	\end{align*}
	Then, a maximum likelihood estimate $\hat{z}^+ \in \mathbb{M}$ of $\Psi(\Flow_{f + \sum_{i=1}^m p_i g_i}^t(x))$ is given by
	\begin{align*}
		\hat{z}^+ \in \pi_{\Sigma^{-1}}(z^+) \qquad\text{with}\qquad
		\Sigma := \sum_{i,j=0}^m Q_{ij} p_i p_j.
	\end{align*}
\end{corollary}
We conclude this section commenting on data-driven approximations of the covariance surrogate $Q$.
\begin{remark}[Data-driven approximation by EDMD]\label{rem:edmdcov}
	To compute the matrix representation of $Q$ serving as an approximation of the covariance matrix~\eqref{eq:koopman_system_noise} similar to EDMD applied to the Koopman operator, we perform a Monte-Carlo approximation of the integrals appearing in \eqref{eq:covariance_map_Q} by a sum over finitely many data points. 
\end{remark}

\subsection{Numerical results}\label{subsec:numres}

In this part, we present various numerical examples to illustrate the proposed parametric Koopman scheme and the corresponding reprojections. We discuss a scalar example in the subsequent Subsection~\ref{subsec:pitchfork}, a two-dimensional Duffing oscillator in Subsection~\ref{subsec:duffing} and the three-dimensional Lorenz dynamics in Subsection~\ref{subsec:lorenz}. In all these examples, we discretize the ODE by means of a fourth-order Runge-Kutta method with adaptive step-size control. For the evaluation, we distinguish three different cases of EDMD-based predictions:
\begin{itemize}
	\item[(i)] Standard prediction without reprojection. %
	\item[(ii)] Prediction with coordinate reprojection as described in Remark~\ref{rem:coord}; applicable if the lifting map is injective.
	\item[(iii)] Maximum likelihood-based closest-point reprojection as proposed and outlined in Subsection~\ref{subsec:reprojections}.
\end{itemize}
In all numerical examples, if not stated otherwise, we approximate the Koopman operator  using a polynomial dictionary up to order five with a sampling time step of $t = 0.1$. %
\new{We recall that, in the sense described in Section~\ref{sec:intro}, the examples investigate whether the surrogate reproduces the bifurcation behaviour for a set of fixed parameter values.}

\subsubsection{Supercritical pitchfork bifurcation}\label{subsec:pitchfork}

First, we consider the scalar system 
\begin{align*}
	\frac{\mathrm{d}}{\mathrm{d}t} x(t) = px(t) - x(t)^3
\end{align*}
with parameter $p\in \mathbb{R}$.
This system corresponds to a supercritical pitchfork bifurcation, as with the parameter becoming positive, two stable equilibria arise. More precisely, we have the following two cases:
\begin{itemize}
	\item[(i)] $p\leq 0$: The only equilibrium is the origin which is asymptotically stable.
	\item[(ii)] $p > 0$: The equilibria are given by $\{0,-\sqrt{p},\sqrt{p}\}$, where the origin is unstable, and $\sqrt{p}$ resp.\ $-\sqrt{p}$ are asymptotically stable for positive resp.\ negative initial values.
\end{itemize}
As a parameter and state domain, we choose $\mathbb{L}=\mathbb{X} = [-2,2]$. %

\begin{figure}[htb]
	\centering
	\includegraphics[width=\linewidth]{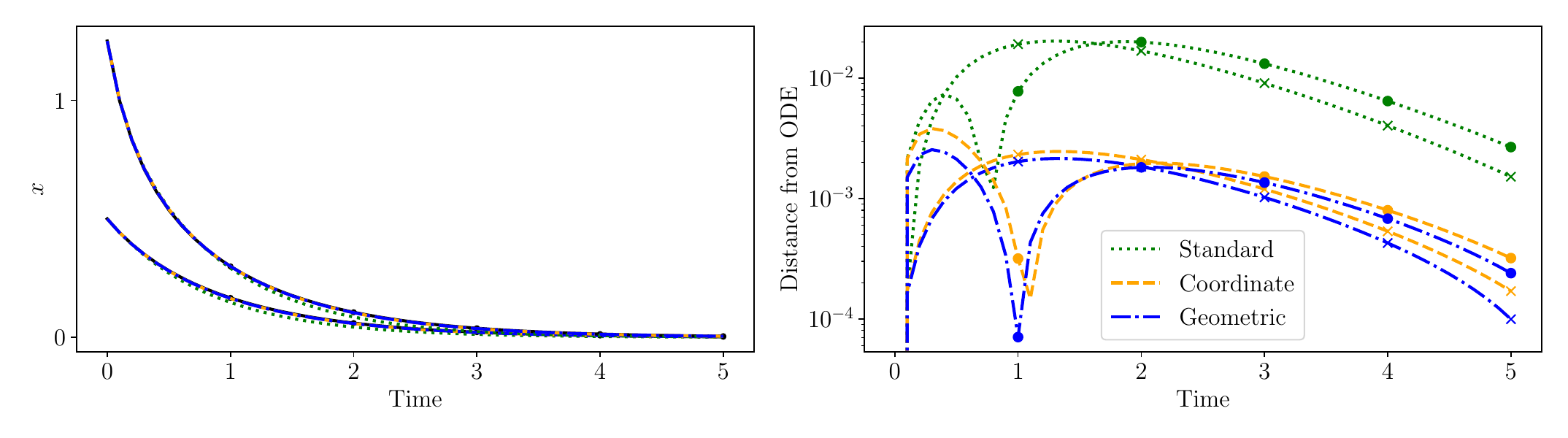}
	\includegraphics[width=\linewidth]{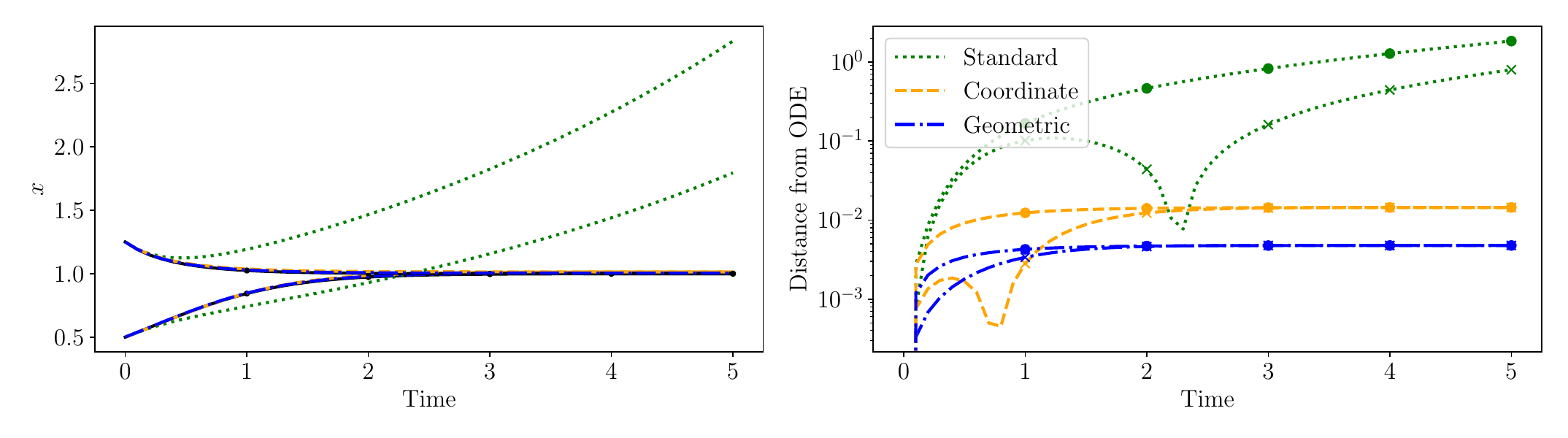}
	\caption{Pitchfork bifurcation for $p=-1$ (top) and $p=1$ (bottom) for initial condition $x^0 = 0.5$ and $x^0 = 1.25$. Left: Trajectories over time. Right: Errors over time with $x^0=0.5$ marked by crosses and $x^0 = 1.25$ marked by bullets.}
	\label{fig:pitchfork}
\end{figure}

In Figure~\ref{fig:pitchfork}, we depict the scalar state (left) and the error over time (right). For $p=-1$, as depicted in the top row, we observe that the globally asymptotically stable nature of the origin renders all three Koopman-based predictions robust, that is, the error is relatively stable over the horizon. 
When transitioning to $p=1$, however, we clearly see that a standard prediction fails after a short time period. In particular, while the formed equilibrium at $x=\sqrt{p}=1$ is correctly identified using reprojections, this is not the case for the standard prediction.

\begin{figure}[htb]
	\centering
	\includegraphics%
	[width=.45\linewidth]
	{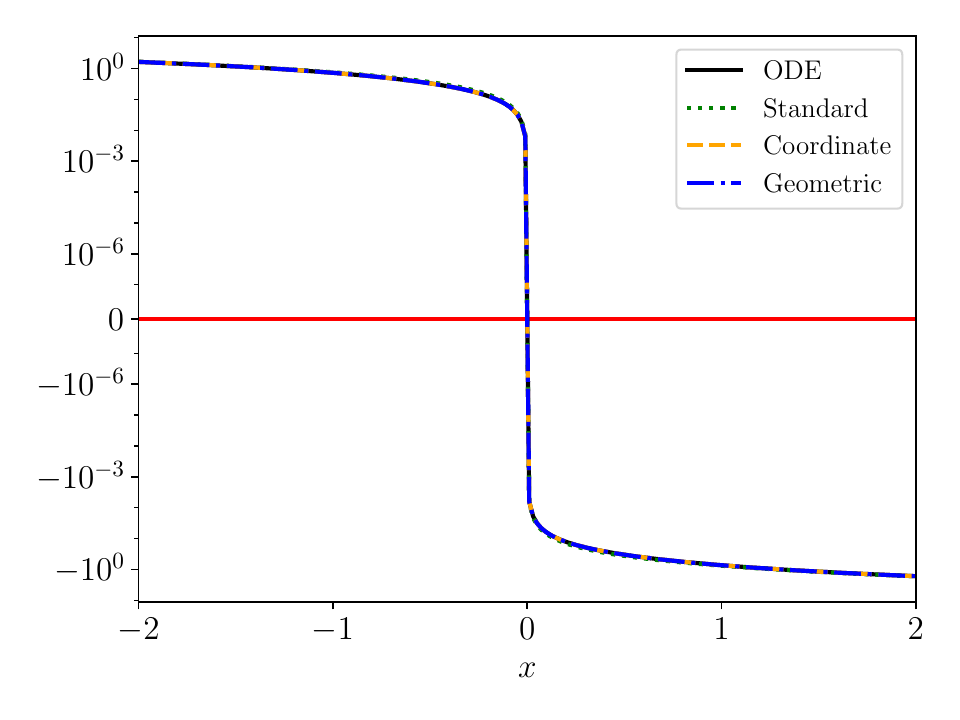}
	\includegraphics%
	[width=.45\linewidth]
	{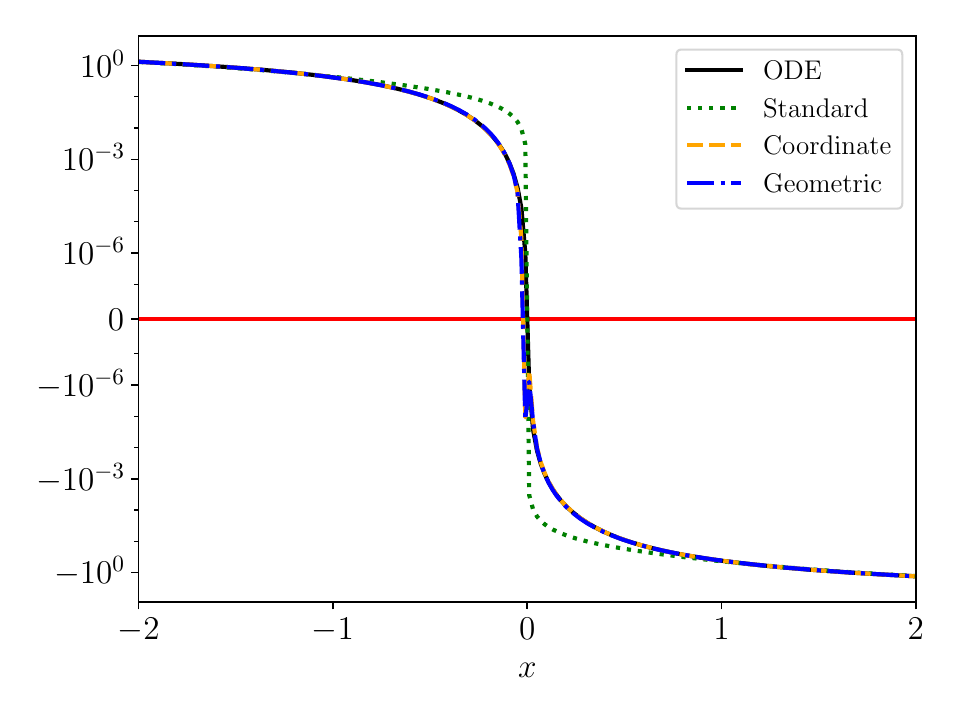}\\
	\includegraphics%
	[width=.45\linewidth]
	{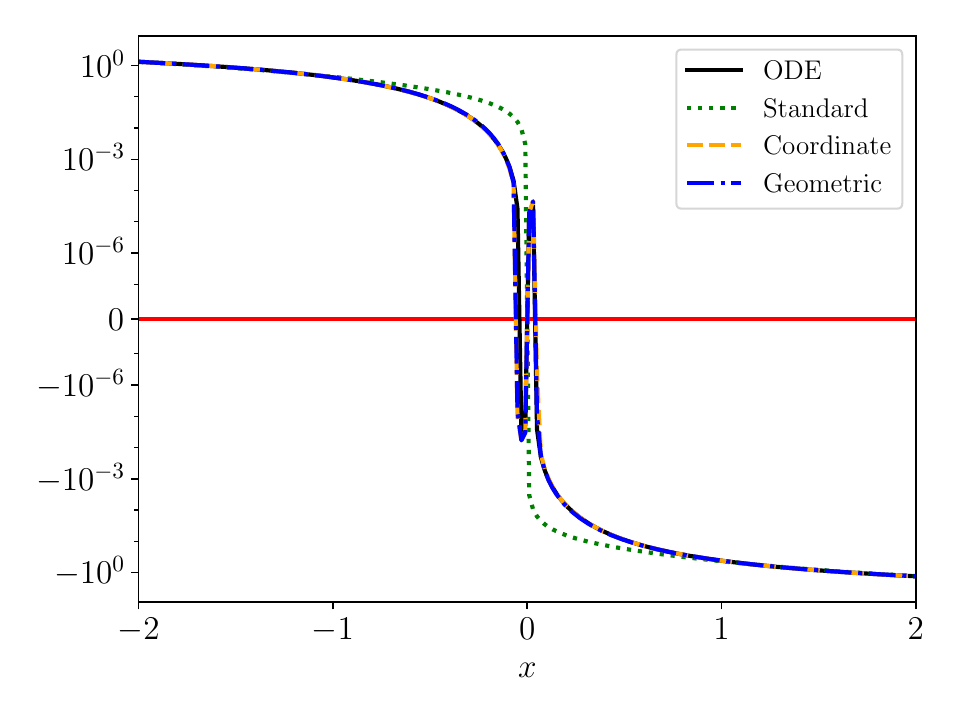}
	\includegraphics%
	[width=.45\linewidth]
	{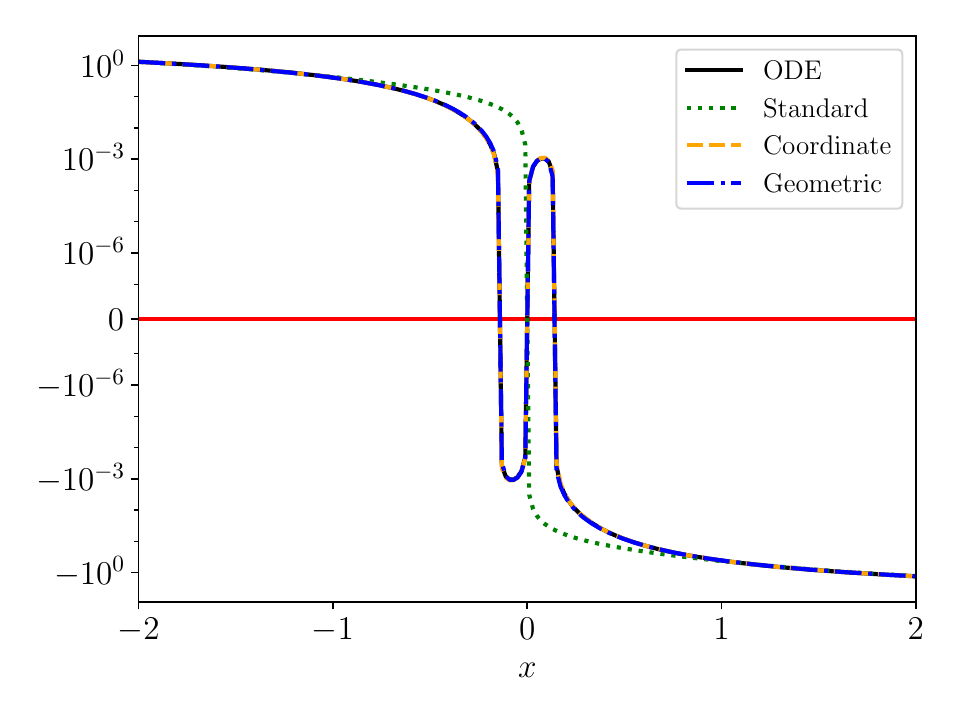}\\
	\includegraphics%
	[width=.45\linewidth]
	{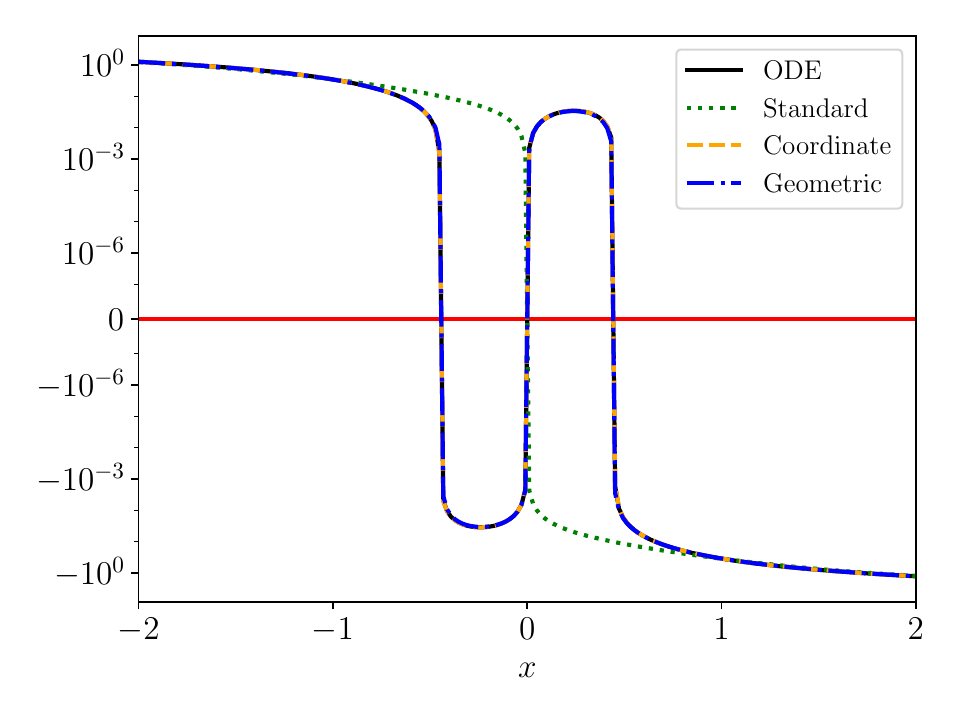}
	\includegraphics%
	[width=.45\linewidth]
	{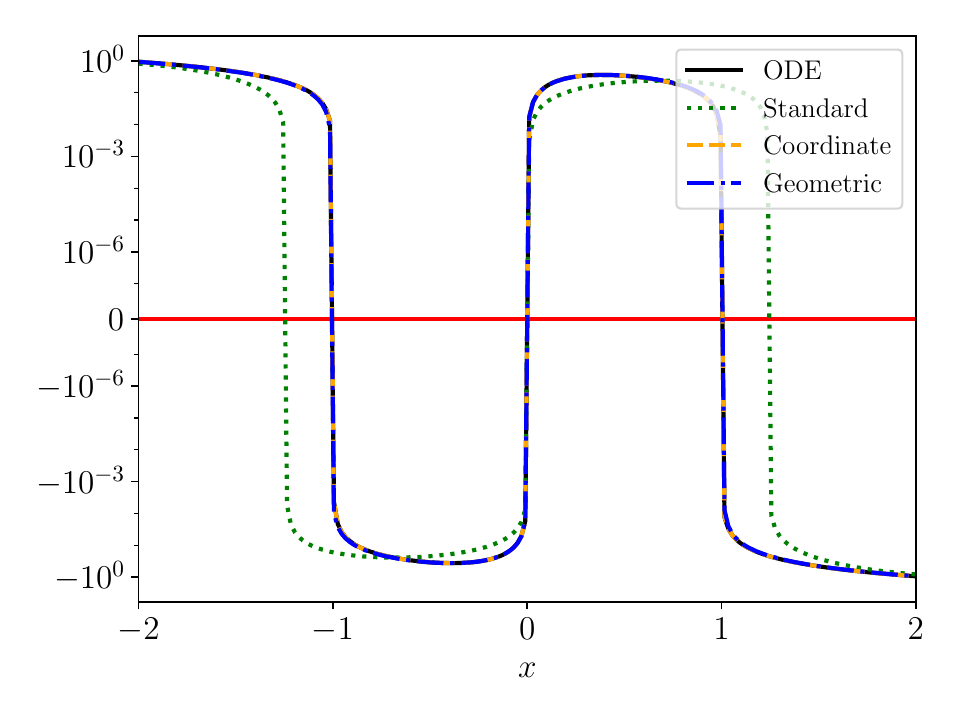}
	\caption{Pitchfork bifurcation: Illustration of bifurcation behaviour for $p\in \{-1,0,0.002,0.02,0.2,1\}$ (row-wise).}
	\label{fig:pitchfork2}
\end{figure}

\noindent In Figure~\ref{fig:pitchfork2}, we illustrate how the bifurcation at $p=0$ is formed. The \new{horizontal red line at zero} illustrates the stability margin: If the dynamics is below this threshold, it is decreasing and if it is above this threshold, it is increasing. Consequently, the fixed points are then the intersection of the dynamics with the \new{horizontal} red line. We observe that, the bifurcation, the emerging fixed points and the corresponding stability behaviour is very well captured by both the coordinate and the maximum likelihood-based reprojection. The standard Koopman prediction does not capture the emerging fixed points and, as can be seen in the bottom left figure, also leads to a completely different stability behaviour in the interval $[-\sqrt{p},\sqrt{p}]$.

\subsubsection{Duffing oscillator}\label{subsec:duffing}
As a second example, we consider the two-dimensional Duffing oscillator
\begin{align*}
	\frac{\mathrm{d}}{\mathrm{d}t}
	\begin{pmatrix}
		x_1(t)\\
		x_2(t)
	\end{pmatrix}
	= \begin{pmatrix}
		x_2(t)\\
		-\delta x_2(t) - \alpha x_1(t) - \beta x_1(t)^3
	\end{pmatrix}            
\end{align*}
with fixed parameters $\delta = 0$ (no friction) and $\beta = 1$. As a varying parameter, we consider the stiffness $\alpha \in \mathbb{L}:=[-2,2]$. 

\begin{figure}[H]
	\centering
	\includegraphics[width=\linewidth]{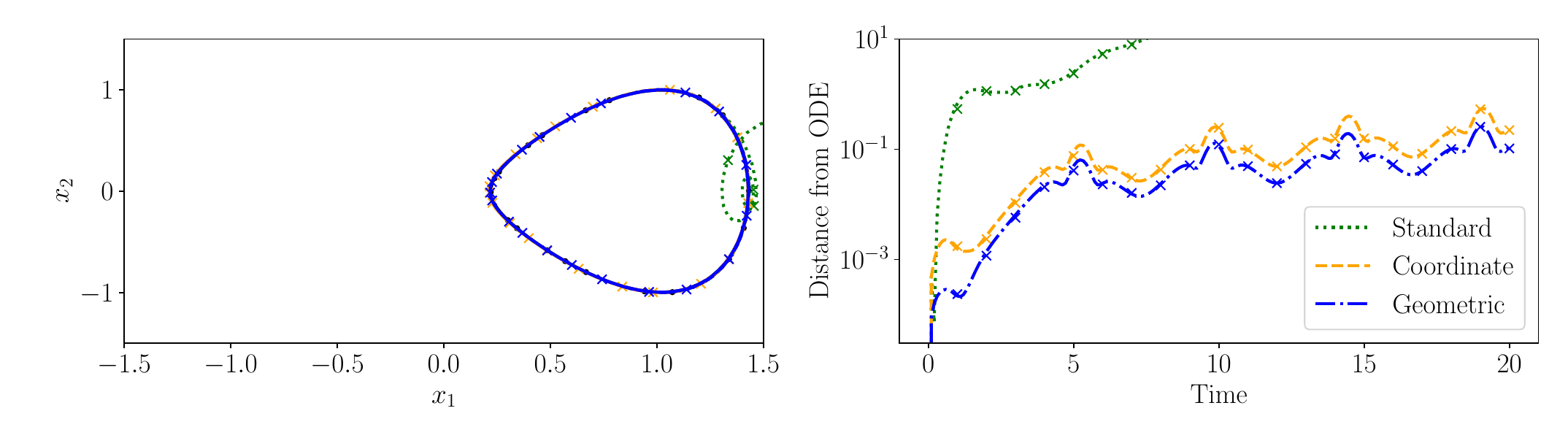}
	\includegraphics[width=\linewidth]{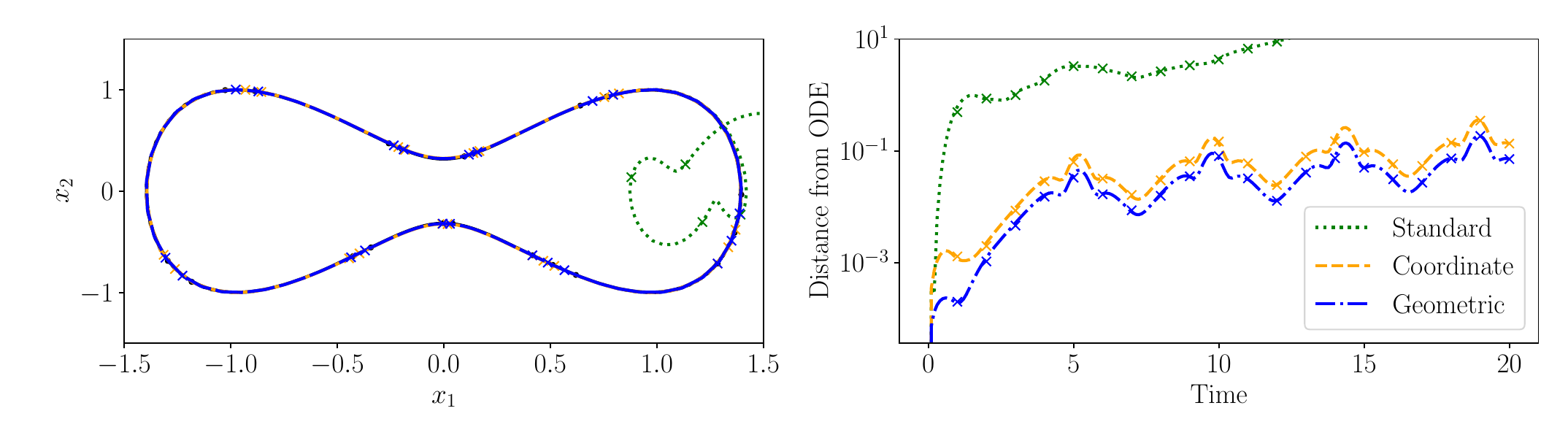}
	\includegraphics[width=\linewidth]{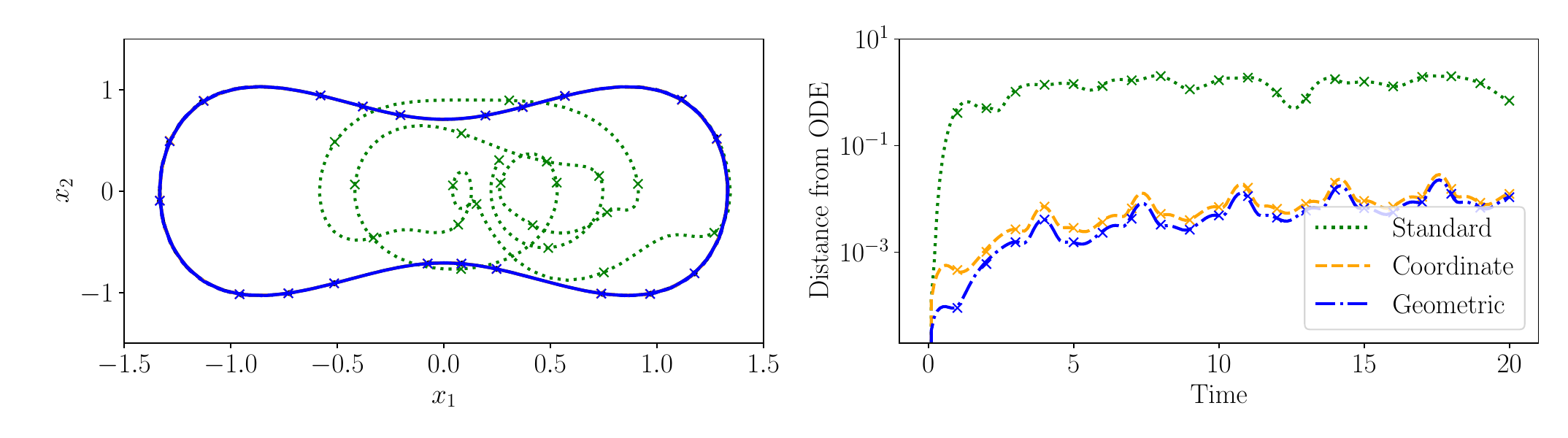}
	\includegraphics[width=\linewidth]{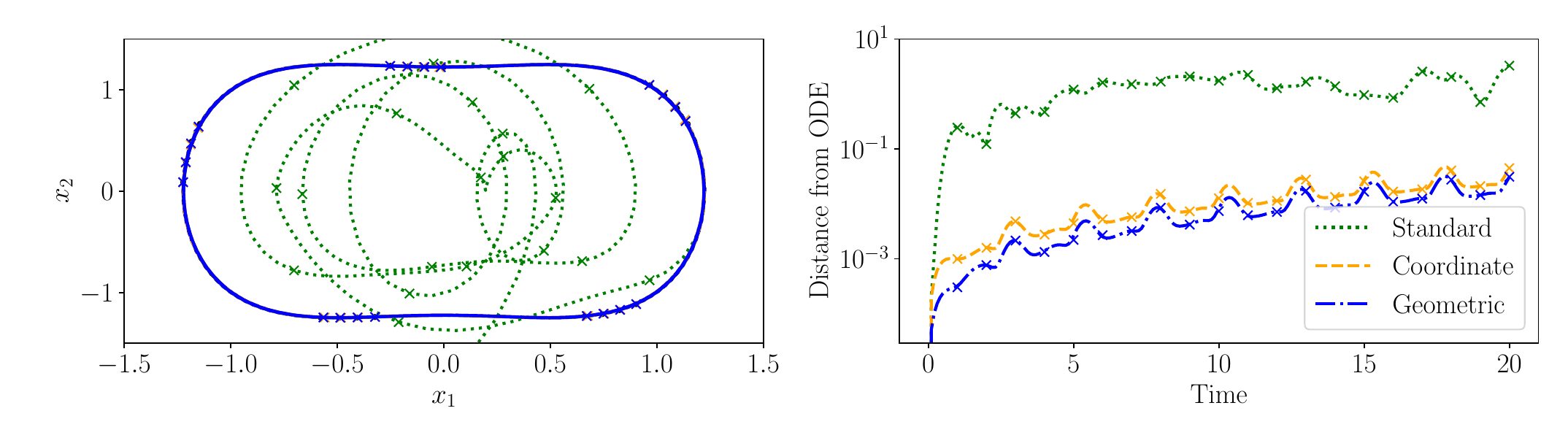}
	\caption{Different predictors for the parametric Duffing oscillator with $\alpha\in \{-2.1,-1.9,-1.5,-0.5\}$ (top to bottom)}
	\label{fig:duffing}
\end{figure}

In Figure~\ref{fig:duffing}, we again plot the performance of the three different Koopman-based predictions for four different stiffness parameters $\alpha$. We observe that the standard prediction fails for all parameter choices, while the projection-based methods both perform well, even for stiffness $\alpha = -2.1$ outside the parameter domain depicted in the top row. In particular, the errors remain relatively stable along the prediction horizon.

\subsubsection{Lorenz attractor}\label{subsec:lorenz}
\begin{figure}[htb]
	\centering
	\includegraphics[width=0.48\linewidth]{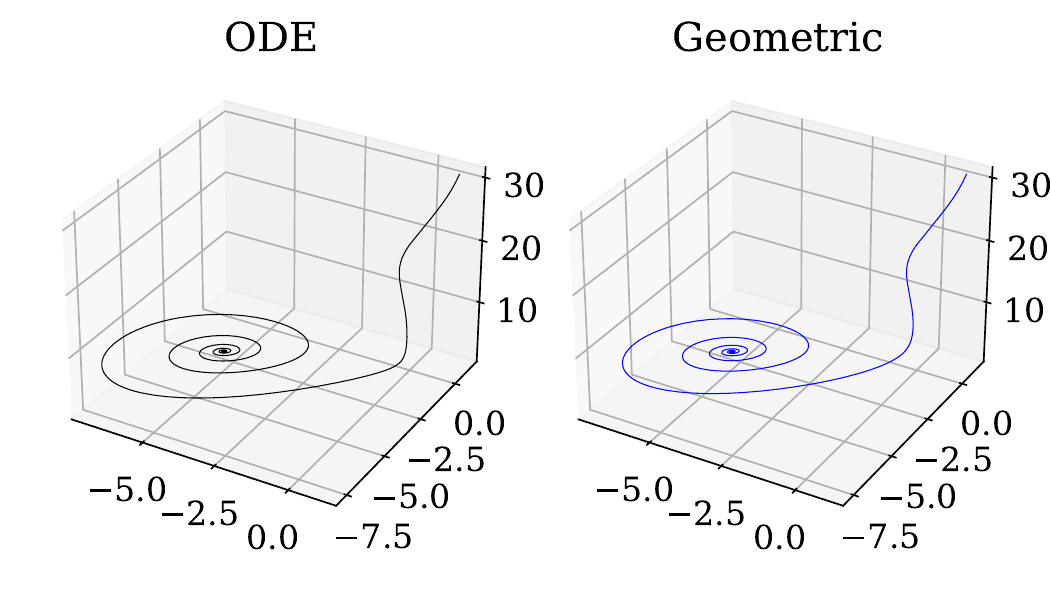}
	\includegraphics[width=0.48\linewidth]{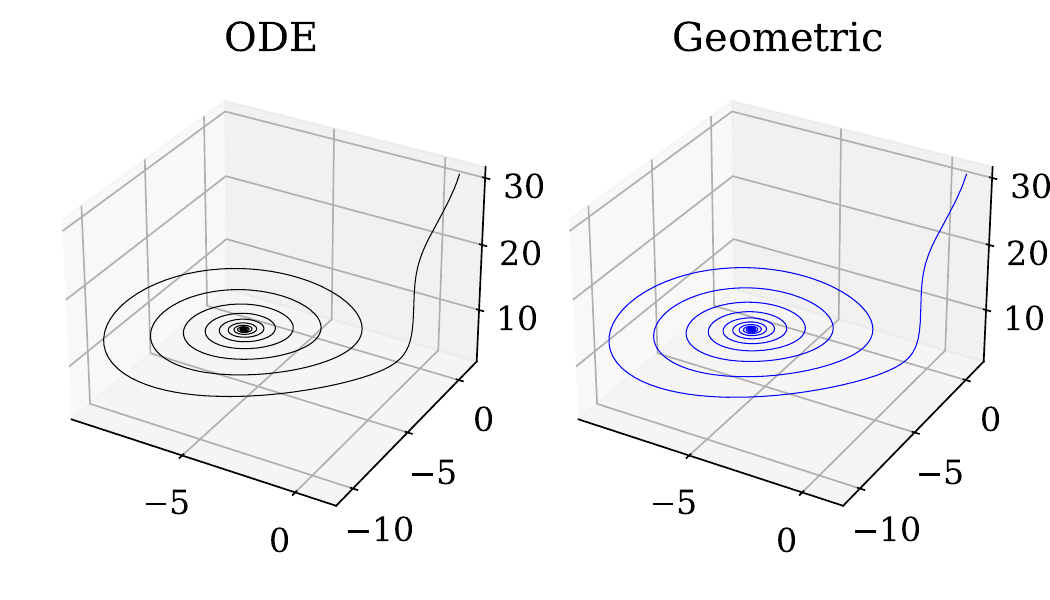}
	\includegraphics[width=0.48\linewidth]{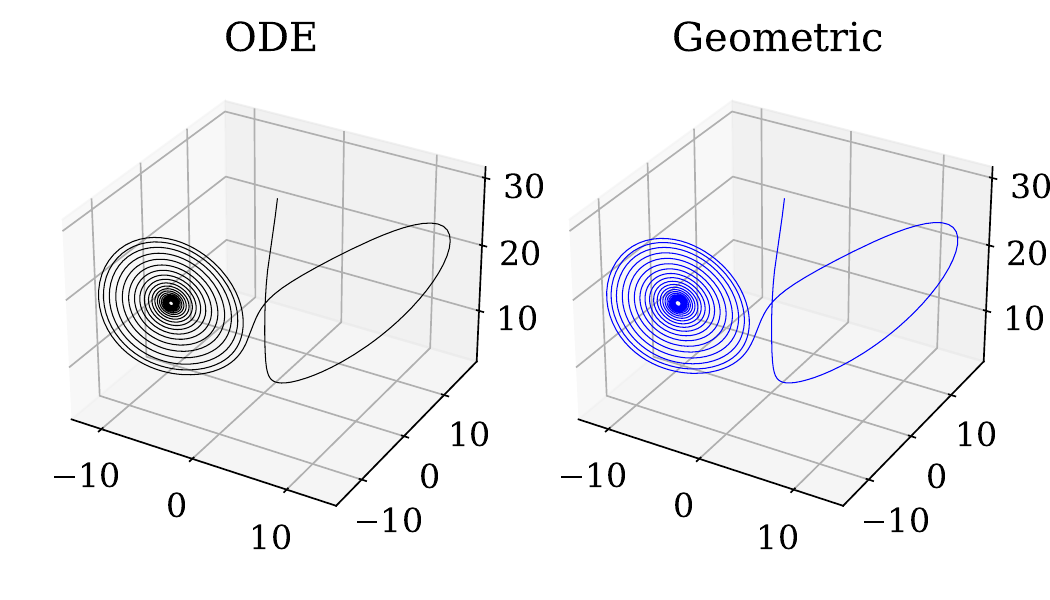}
	\includegraphics[width=0.48\linewidth]{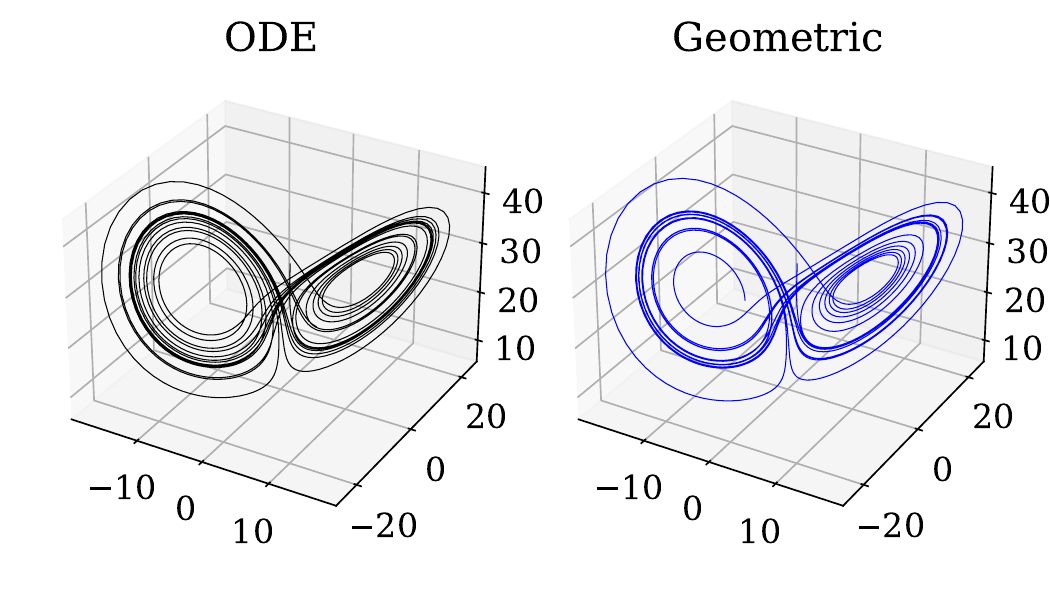}
	\caption{Formation of the Lorenz attractor by variations of the parameter $\rho \in \{8,12,18,28\}$ (row-wise).}
	\label{fig:lorenz}
\end{figure}

As a last example, we consider the three-dimensional Lorenz system given by
\begin{align}\nonumber %
	\dot{x}_1 &= \sigma (x_2-x_1), &
	\dot{x}_2 &= x_1(\rho-x_3)-x_2, &
	\dot{x}_3 &= x_1x_2 - \beta x_3
\end{align}
on the domain $\mathbb{X} = [-20,20]\times[-20,20]\times[10,50]$. We fix $\sigma = 10$ and $\beta = 8/3$ and consider $\rho\in [10,30]$ as a parameter. Note that the most common choice of the parameter giving rise to the well-known chaotic attractor is given by $\rho = 28$. As a time step, we choose $t=0.01$ to capture the fast dynamics of the system and as a dictionary, we choose a monomial dictionary up to order three. In Figure~\ref{fig:lorenz}, we illustrate the formation of the attractor for different choices of the parameter $\rho$ for the geometric reprojection. While the propagation using the standard Koopman model diverges already after a few time steps, we observe that the attractor is very well captured by the proposed geometric reprojection. The coordinate reprojection behaves very similarly to the geometric reprojection, which is why we only show the latter in Figure~\ref{fig:lorenz}.

\begin{figure}[htb]
	\centering
	\includegraphics[width=0.48\linewidth]{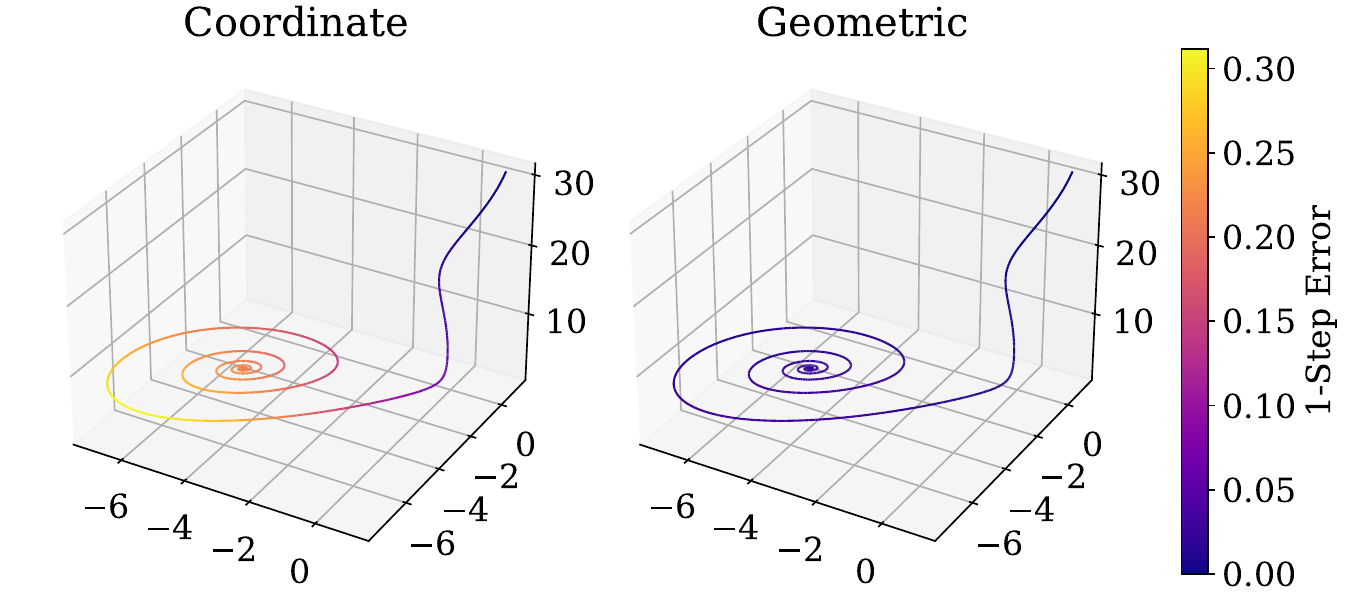}
	\includegraphics[width=0.48\linewidth]{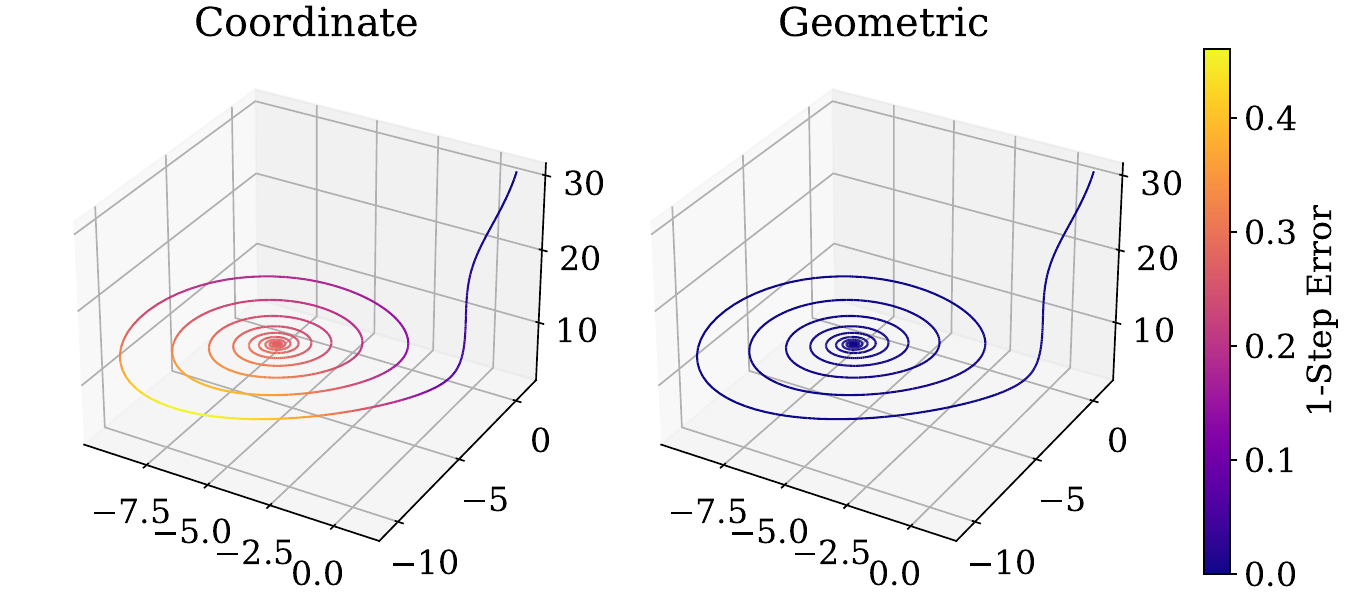}
	\includegraphics[width=0.48\linewidth]{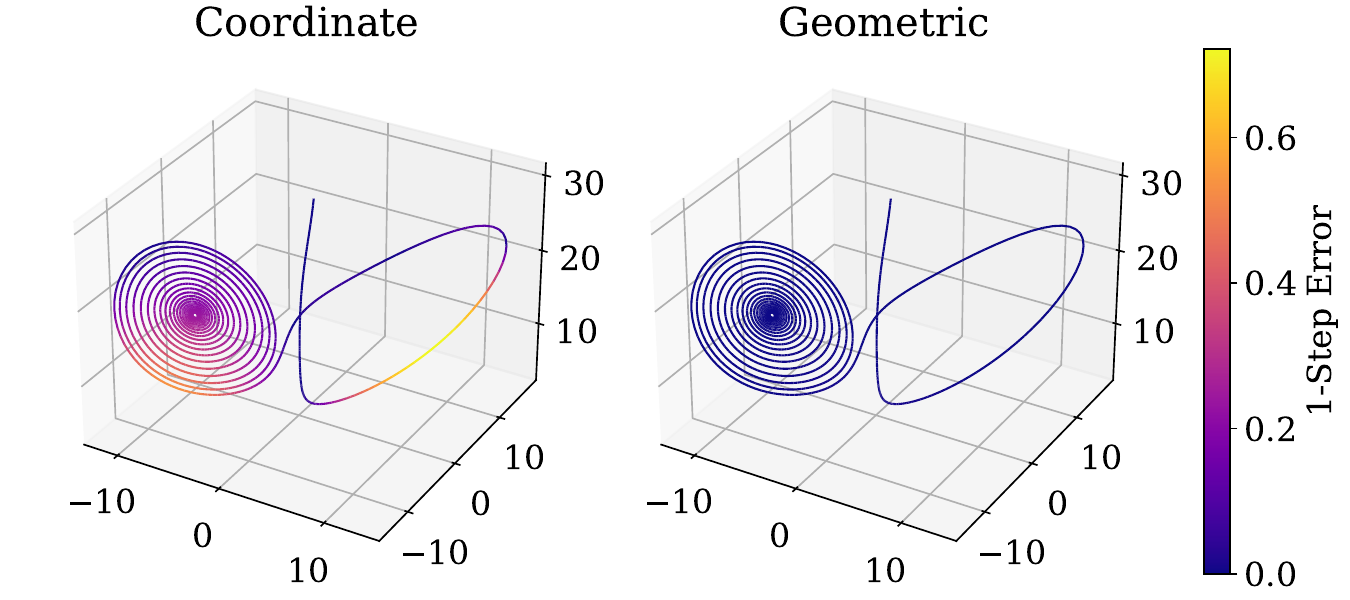}
	\includegraphics[width=0.48\linewidth]{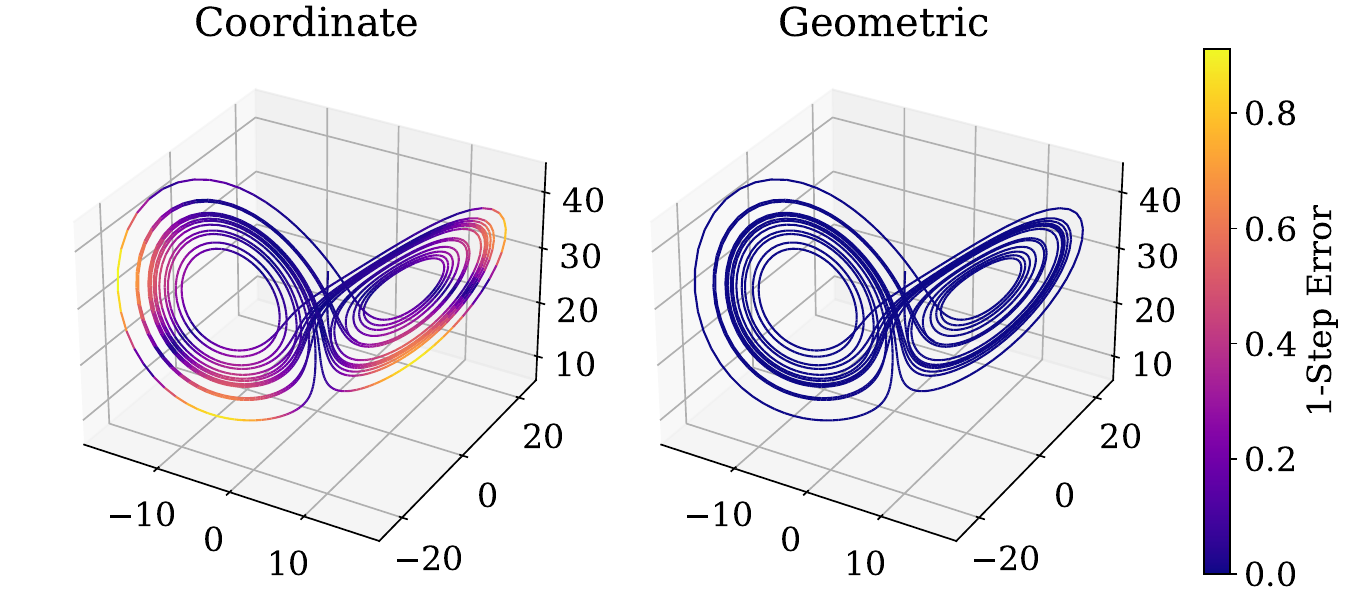}
	\caption{Formation of the Lorenz attractor by variations of the parameter $\rho \in \{8,12,18,28\}$ (row-wise). Monomial dictionary without coordinate map of the first state, i.e., $\mathbb{V}_3 \setminus \{x_1\}$. }
	\label{fig:lorenz_nocoord}
\end{figure}

In Figure~\ref{fig:lorenz_nocoord}, we illustrate the results using a monomial dictionary without the first coordinate map. To obtain state predictions (and to define a coordinate projection), one may leverage injectivity of the lifting. More precisely, in the spirit of Remark~\ref{rem:coord}, we may define a weighting matrix $W_c\in \R^{M\times M}$ by setting the diagonal entries corresponding to $x_2,x_3$ and $x_1^3$ to one, and the remainder to zero. We observe in Figure~\ref{fig:lorenz_nocoord} that the geometric reprojection outperforms this coordinate reprojection. While long-term predictions depicted in Figure~\ref{fig:lorenz} worked well for both the coordinate and the geometric reprojections in case of all coordinate functions contained in the dictionary $\mathbb{V} = \mathbb{V}_3$, in case of $\mathbb{V}=\mathbb{V}_3\setminus {x_1}$, the larger error of the coordinate reprojection (depicted in Figure~\ref{fig:lorenz_nocoord}) leads to divergence already after a few time steps. Hence, when considering a dictionary that does not contain the coordinate functions, the geometric projection is the method of choice.

\section{Efficient implementation}\label{sec:efficient}
\noindent 
While the previous investigations showcased a favourable influence of the reprojection on the stability and accuracy of Koopman-based predictions, each evaluation requires the solution of a nonlinear minimization problem. In this section, we show how the numerical effort may be kept tractable by (i) reducing the number of projections required and (ii) providing an efficient numerical method to compute the projections.

To address (i), we propose an adaptive choice of the reprojection time in Subsection~\ref{subsec:multistep}. In particular, we will follow the linear model for several time steps and then evaluate the projection only after a few time steps by means of an error estimator. To address (ii), we formulate a Riemannian Newton method in Subsection~\ref{subsec:newton}. The nature of the problems to be solved being a series of optimization problems along a trajectory provides a very natural candidate for an initial guess enabling fast local convergence of the Newton method.

\subsection{Multistep predictions}\label{subsec:multistep}

In this part, we employ the covariance model to provide a triggering condition for the reprojection step. This allows us to conduct a reprojection only if needed, similarly to event-triggered mechanisms for feedback linearization~\cite{umlauft2019feedback}, or learning in general~\cite{solowjow2020event} conceptually following the sequential exploration-exploitation framework~\cite{lew2022safe}.

When considering data-driven multistep predictions, the lifted state $\hat{z}(k)$ follows
\begin{align*}
	\hat{z}(k+1) = \left(\widehat{K}_0 + \sum_{i=1}^m p_i \widehat{K}_i \right) \hat{z}(k) \qquad k=0,1,\ldots
\end{align*}
where $\hat{z}(0) = \Psi(x(0))\in \mathbb{M}$ for an initial state $x(0)$.
From the approximation \eqref{eq:koopman_system_noise}, we have that
\begin{align*}
	\left(\widehat{K}_0 + \sum_{i=1}^m p_i \widehat{K}_i\right) \Psi(x(k)) \approx \Psi(x(k+1)) + \mu,
\end{align*}
where $\mu \sim N\left(0, \sum_{i,j=0}^m Q_{ij} p_i p_j\right)$. %

Let $k\in \mathbb{N}$ be fixed and consider a prediction $\hat{z}(k)\in \R^M$ such that $ \hat{z}(k) \approx \Psi(x(k)) + \eta_k$, where $\eta_k \sim N(0, \Sigma_k)$ for a given covariance $\Sigma_k$.
Then,
\begin{align*}
	\hat{z}(k+1)
	&= \left(\widehat{K}_0 + \sum_{i=1}^m p_i \widehat{K}_i \right) \hat{z}(k) \\
	&= \left(\widehat{K}_0 + \sum_{i=1}^m p_i \widehat{K}_i \right) (\Psi(x(k)) + \eta_k) \\
	&\approx \Psi(x(k+1)) + \new{\mu} + \left(\widehat{K}_0 + \sum_{i=1}^m p_i \widehat{K}_i \right) \eta_k.
\end{align*}
It follows that 
\begin{align}
	\hat{z}(k+1) &\approx \Psi(x(k+1)) + \eta_{k+1}\nonumber \\
	\eta_{k+1} &\sim N \left(0, \; \sum_{i,j=0}^m Q_{ij} p_i p_j + \left(\widehat{K}_0 + \sum_{i=1}^m p_i \widehat{K}_i \right) \Sigma_k \left(\widehat{K}_0 + \sum_{i=1}^m p_i \widehat{K}_i \right)^\top \right).\label{eq:multicov}
\end{align}
In other words, the statistics of the error term $\eta_{k+1}$ can be propagated using the parameters and the surrogate models. One can choose to perform reprojection only as needed, for example, when the covariance of $\eta_k$ crosses a threshold. Note that to this end, we have to define a measure of the covariance, e.g., the trace of the matrix or diagonal entries corresponding to the covariance of particular variables.

\begin{figure}[htb]
	\centering
	\includegraphics[width=.9\linewidth]{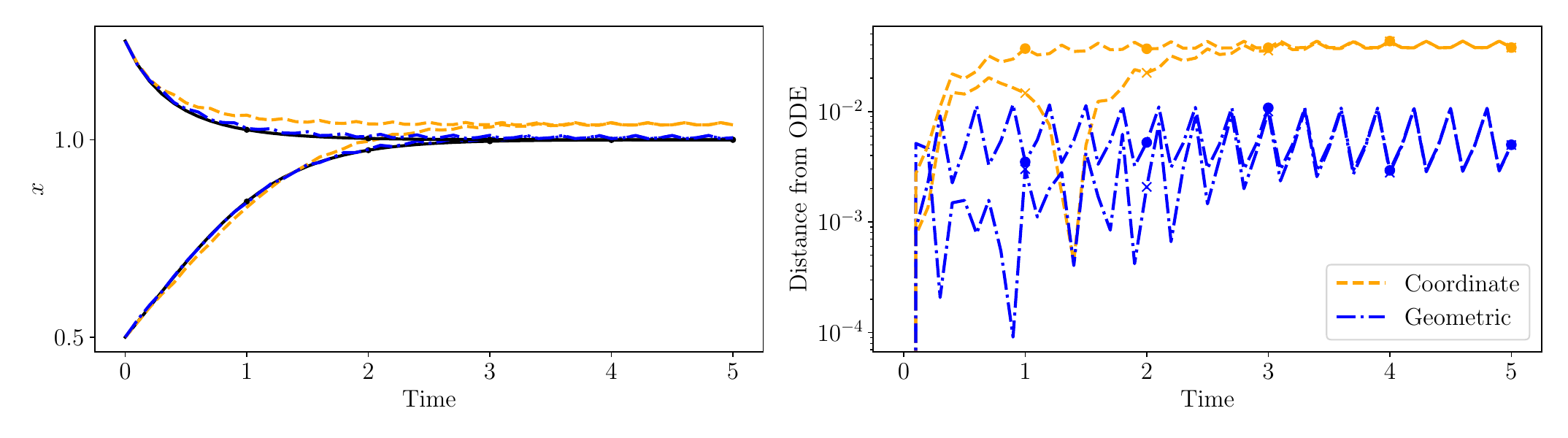}
	\includegraphics[width=.9\linewidth]{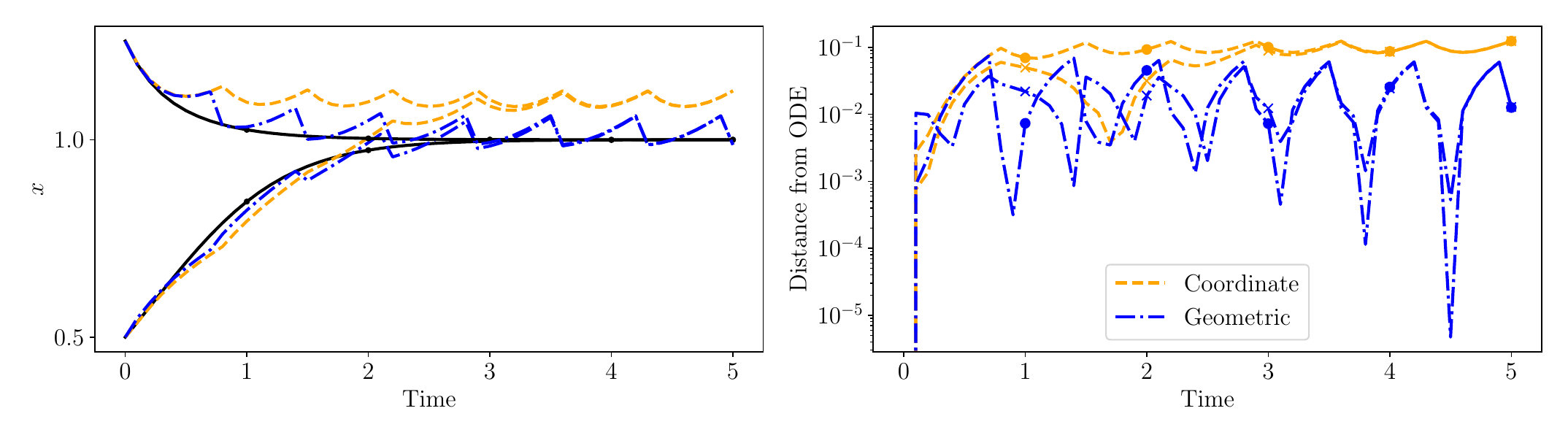}\\
	\includegraphics[width=.9\linewidth]{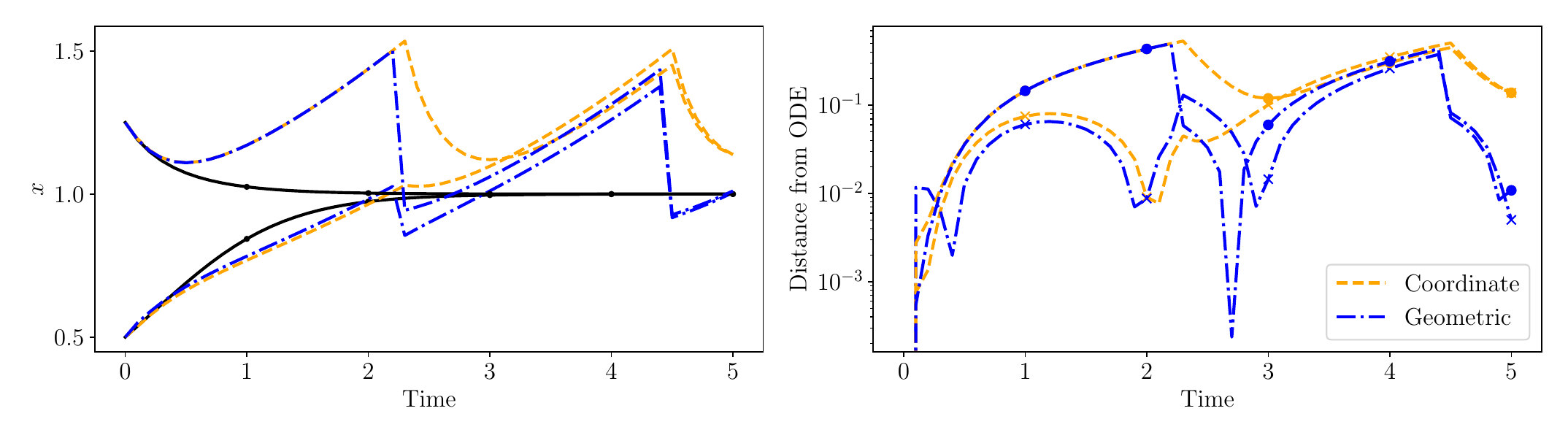}
	\caption{Pitchfork bifurcation with $p=1$: Multistep predictions for a tolerance for the increase of the covariance by a factor in $\{10,100,1000\}$ leading to multistepping with $\{2,6,21\}$ steps (top to bottom). The error on the right is depicted for $x^0=0.5$ by crosses and for $x^0=1.25$ by bullets.}
	\label{fig:multistep_1d}
\end{figure}

In Figure~\ref{fig:multistep_1d}, we depict the performance of the proposed multistep scheme for the scalar example of Subsection~\ref{subsec:pitchfork}. To determine the prediction length, we only check the increase of the $(1,1)$-entry of the covariance matrix corresponding to the coordinate function.

\begin{figure}[htb]
	\centering
	\includegraphics[width=.9\linewidth]{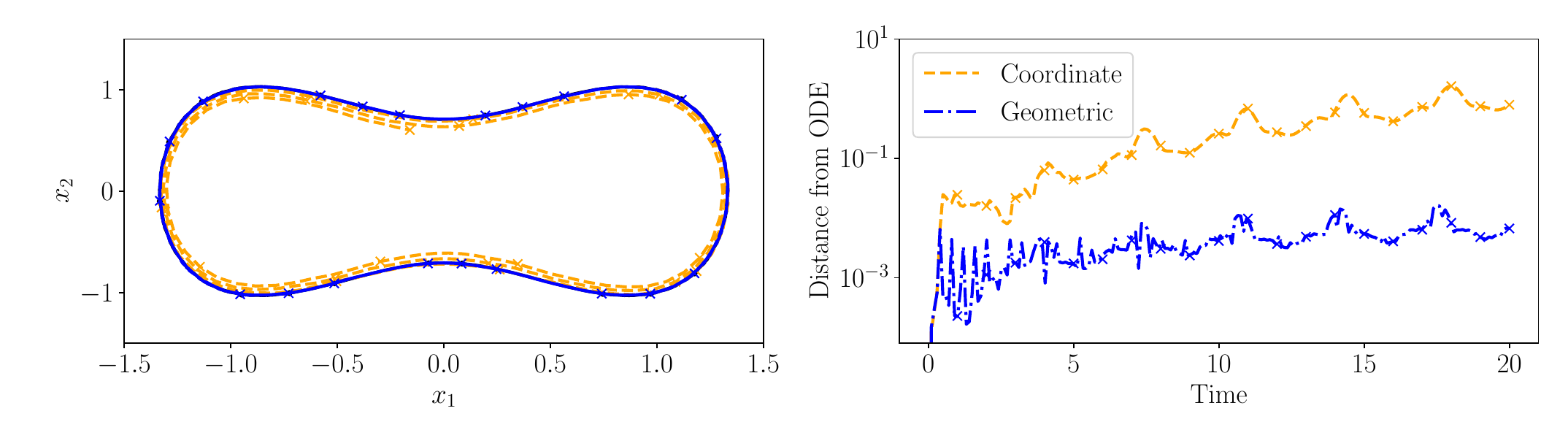}
	\includegraphics[width=.9\linewidth]{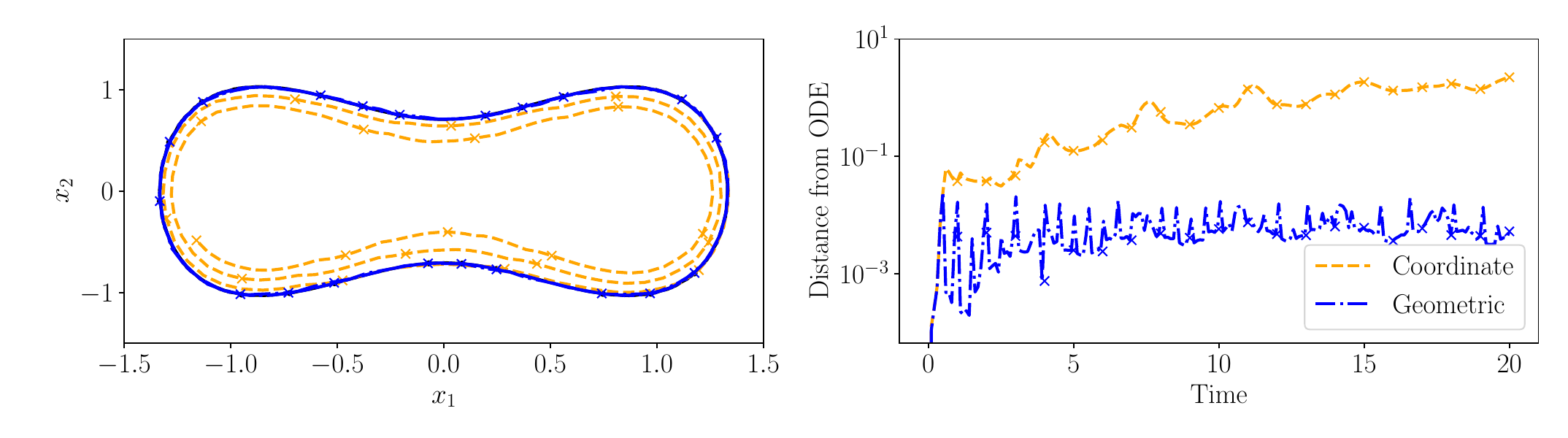}\\
	\includegraphics[width=.9\linewidth]{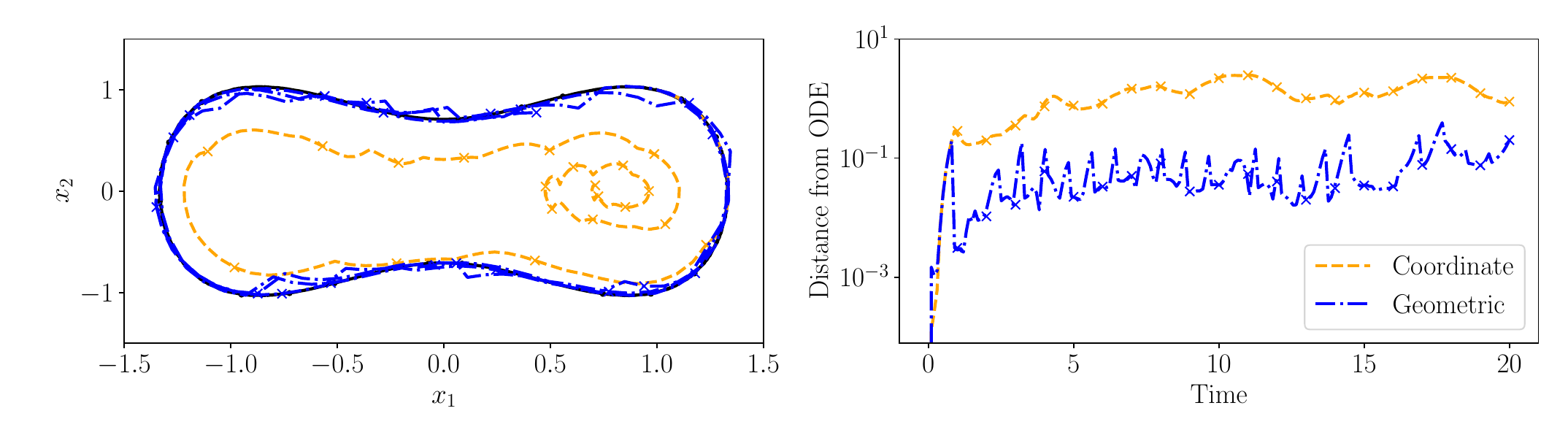}
	\caption{Duffing oscillator with $\alpha = -1.5$: Multistep predictions for a tolerance for the increase of the covariance by a factor in $\{10,20,30\}$ leading to multistepping with $\{3,4,7\}$ steps (top to bottom).}
	\label{fig:multistep_2d}
\end{figure}

In Figure~\ref{fig:multistep_2d}, we illustrate the multistep predictions for the two-dimensional Duffing oscillator of Subsection~\ref{subsec:duffing}. Here, we evaluate the trace of the covariance matrix and reproject if this uncertainty measure has increased by a pre-defined factor.

\new{We briefly summarize the findings of this subsection. The covariance recursion~\eqref{eq:multicov} provides an error indicator that is available online at negligible cost, as it only involves the surrogate matrices and the bilinear form~$Q$, both of which are computed offline. Hence, accuracy may be traded against the number of reprojections in a controlled manner. More precisely, for the pitchfork example depicted in Figure~\ref{fig:multistep_1d}, relaxing the admissible increase of the covariance from a factor of $10$ to a factor of $1000$ increases the number of prediction steps between two reprojections from $2$ to $21$, i.e., it reduces the number of reprojection problems to be solved by roughly an order of magnitude, at the price of an increase of the maximal prediction error by about one to two orders of magnitude. In all cases, the qualitative behaviour, i.e., convergence towards the correct equilibrium, is preserved. For the two-dimensional Duffing oscillator depicted in Figure~\ref{fig:multistep_2d}, the geometric reprojection reproduces the limit cycle even for the largest tolerance considered, whereas the coordinate reprojection fails qualitatively and collapses onto a spurious equilibrium. Hence, the admissible number of steps between two reprojections depends on the problem at hand, and the covariance-based indicator rather than a fixed reprojection frequency is the appropriate quantity to base this decision on.}

\subsection{Efficient Reprojection via Riemannian Newton Method}\label{subsec:newton}

In general, it is not straightforward to implement the proposed reprojection function~\eqref{eq:reprojection_function}. 
To use %
{derivative-based} methods for solving the corresponding \emph{manifold-constrained nonlinear optimization problem}, we suppose %
that the image of~$\Psi$ defined in \eqref{eq:M} has the structure of an embedded manifold.

\begin{definition}[Embedded Submanifold~{\cite[Chapter 4,5]{Lee12}}]
	An \emph{embedding} of a manifold $\mathbb{M}$ into $\R^M$ is a smooth immersion $F : \mathbb{M} \to \R^M$ that is also a topological embedding.
	That is, $F$ is smooth, the derivative $\tD F(x)$ has full rank for all $x \in M$, and $F$ is a homeomorphism onto its image $F(M) \subseteq \R^M$.
	A manifold $\mathbb{M}$ is called an \emph{embedded submanifold (of $\R^M$)} if it is a subset $\mathbb{M} \subseteq \R^M$ and the inclusion map is an embedding.
\end{definition}

\begin{assumption}
	The map $\Psi : \mathbb{X} \to \R^M$ is an embedding; equivalently, $
	\mathbb{M}= \Psi(\mathbb{X}) \subset \R^M$ is an embedded submanifold.
\end{assumption}

To apply optimization methods on the manifold $
\mathbb{M}$, both a Riemannian metric on $\mathbb{M}$ and a retraction from the tangent bundle $\mathrm{T} 
\mathbb{M}$ are required.
We recall the definitions as they are needed.

\begin{definition}[Riemannian Metric~{\cite[Chapter 13]{Lee12}}]
	For a manifold $\mathbb{M}$, a \emph{Riemannian metric on $\mathbb{M}$} is a 2-
	tensor field $g$ that is positive definite at each point.
	That is, for each $p \in \mathbb{M}$, the metric defines a positive definite bilinear map $g_p : \tT_p \mathbb{M} \times \tT_p \mathbb{M} \to \R$.
\end{definition}

The Riemannian metric may be simply inherited from the ambient Euclidean space; that is, for any $p\in \mathbb{M}$ and $\eta_1, \eta_2 \in \tT_p \mathbb{M}$, the Riemannian metric at $p$ is defined by
\begin{align}\label{eq:riemannian_metric}
	g_p(\eta_1, \eta_2)  := \langle \eta_1, \eta_2 \rangle,
\end{align}
where the right-hand side is the ordinary inner product in $\R^M$ applied by the inclusion of $\tT_p \mathbb{M}$ as a subspace of $\R^M$.

\begin{definition}[Retraction~{\cite[Definition 4.1.1]{AbsiMaho09}}]
	A \emph{retraction} on a manifold $\mathbb{M}$ is a smooth map $R: \tT \mathbb{M} \to \mathbb{M}$ satisfying the following two properties, where $R_p : \tT_p \mathbb{M} \to \mathbb{M}$ is the restriction of $R$ to $\tT_p \mathbb{M}$:
	\begin{enumerate}
		\item $R_p(0_p) = p$, where $0_p \in \tT_p \mathbb{M}$ is the zero vector in the tangent space at $p$.
		\item Through the canonical identification $\tT_{0_p} \tT_p \mathbb{M} \simeq \tT_p \mathbb{M}$, the differential of $R_p$ satisfies
		\begin{align*}
			\tD R_p(0_p) = \mathrm{id}_{\tT_p \mathbb{M}},
		\end{align*}
		where $\mathrm{id}_{\tT_p \mathbb{M}} : \tT_p \mathbb{M} \to \tT_p \mathbb{M}$ denotes the identity map.
	\end{enumerate}
\end{definition}

The definition of $\mathbb{M}$ as the image of an embedding $\Psi$ lends itself to defining a retraction.
\begin{lemma}\label{lem:retraction}
	Consider the manifold $\mathbb{M} = \Psi(\mathbb{X}) \subset \R^M$.
	Define the map $R : \tT \mathbb{M} \to \mathbb{M}$ by
	\begin{align*}
		R_p(\eta) = \Psi(\Psi^{-1}(p) + \tD \Psi(\Psi^{-1}(p))^{-1} \eta),
	\end{align*}
	for all $p \in \mathbb{M}$ and $\eta \in \tT_p \mathbb{M}$.
	Then $R$ is a retraction.
\end{lemma}

\begin{proof}
	Smoothness of $R$ is immediate as it is a composition of smooth maps.
	Let $p \in \mathbb{M}$ be arbitrary and denote the zero element of the tangent space as $0_p \in \tT_p \mathbb{M}$.
	Then the retraction of the zero element satisfies
	\begin{align*}
		R_p(0_p) = \Psi(\Psi^{-1}(p) + \tD \Psi(\Psi^{-1}(p))^{-1} 0_p)
		= \Psi(\Psi^{-1}(p) + 0) = p.
	\end{align*}
	Finally, for any $\eta \in \tT_{0_p} \tT_p \mathbb{M} \simeq \tT_p \mathbb{M}$,
	\begin{align*}
		\tD R_p(0_p)[\eta]
		&= \at{\ddt}{t=0} R_p(0_p + t \eta) \\
		&= \at{\ddt}{t=0} \Psi(\Psi^{-1}(p) + \tD \Psi(\Psi^{-1}(p))^{-1} (0_p + t \eta)) \\
		&= \at{\ddt}{t=0} \Psi(\Psi^{-1}(p) + t \tD \Psi(\Psi^{-1}(p))^{-1} \eta) \\
		&= \tD \Psi(\Psi^{-1}(p))\tD \Psi(\Psi^{-1}(p))^{-1} \eta \\
		&= \eta.
	\end{align*}
	Thus $\tD R_p(0_p)$ is exactly the identity map on $\tT_p \mathbb{M}$ for all $p$.
	This shows that $R$ is indeed a retraction.
\end{proof}

Let $R : \tT \mathbb{M} \to \mathbb{M}$ be the retraction defined in Lemma \ref{lem:retraction}.
Then we may apply a Riemannian Newton method \cite[Algorithm 5]{AbsiMaho09} to find a local solution to the reprojection map $\pi_W$.
Specifically, given a point $z \in \R^M$, define
\begin{align*}
	l_z : \mathbb{M} \to \R, && \new{l_z}(p) := \frac{1}{2} \vert p - z \vert_W^2.
\end{align*}
Then the Newton equation is
\begin{align}\label{eq:newton_equation}
	\Hess l_z(p) \eta = - \grad l_z(p).
\end{align}
The gradient and Hessian of $l_z$ are easily computed to be
\begin{align*}
	\grad l_z(p) &= W(p-z) &
	\Hess l_z(p) &= W,
\end{align*}
using the embedded representation of $\tT_p \mathbb{M}$.
The Newton equation \eqref{eq:newton_equation} thus appears trivial to solve by $\eta = (p-z)$.
However, this neglects the fact that the Hessian and gradient of $l_z$ are restricted to the tangent space of $\mathbb{M}$ at $p$.
This issue can be addressed by finding a minimal basis of the tangent space at $p$.

Let $x = \Psi^{-1}(p)$.
Then the tangent space $\tT_p \mathbb{M}$ is exactly the image of the differential $\tD \Psi(x) : \R^d \to \R^M$, and this map is injective.
Therefore any $\eta \in \tT_p \mathbb{M}$ may be written as $\tD \Psi(x) v$ for a unique $v \in \R^d$, and the standard basis vectors $e_1,\ldots,e_d$ of $\R^d$ induce a basis $\tD \Psi(x) e_1, \ldots \tD \Psi(x) e_d \in \R^M$ for $\tT_p \mathbb{M}$.
In this basis, the gradient and Hessian may be written as
\begin{align*}
	\grad l_z(p) &= \tD \Psi(x)^\top W(p-z) &
	\Hess l_z(p) &= \tD \Psi(x)^\top W \tD \Psi(x).
\end{align*}
Let $v \in \R^d$ denote the unique vector in $\R^d$ such that $\eta = \tD \Psi(x) v$. 
Then the Newton equation \eqref{eq:newton_equation} simplifies to
\begin{align}\label{eq:simplified_newton_equation}
	\tD \Psi(x)^\top W \tD \Psi(x) v &= - \tD \Psi(x)^\top W(p-z), \\
	v &= - (\tD \Psi(x)^\top W \tD \Psi(x))^\dag \tD \Psi(x)^\top W(p-z).
\end{align}
The retraction of $\eta$ can also then be simplified to
\begin{align*}
	R_p(\eta)
	&= \Psi(\Psi^{-1}(p) + \tD \Psi(\Psi^{-1}(p))^{-1} \eta), \\
	&= \Psi(\Psi^{-1}(p) + \tD \Psi(\Psi^{-1}(p))^{-1} \tD\Psi(x) v), \\
	&= \Psi(\Psi^{-1}(p) + v).
\end{align*}
Finally, this leads to Algorithm \ref{alg:newton_method} for obtaining $\pi_W(z)$.

\begin{algorithm}
	\caption{Riemannian Newton method for reprojection \eqref{eq:reprojection_function}}
	\label{alg:newton_method}
	\begin{algorithmic}[1]
		\Require $W \in \R^{M\times M}$ positive semidefinite; $z \in \R^M$ the target point.
		\Goal Find $p \in \mathbb{M}$ to minimize $\vert p - z \vert^2_W$.
		\Input An initial guess $p_0 \in \mathbb{M}$.
		\Output An estimate $\hat{p}$ of $\pi_W(z)$.
		\State Set $x_0 = \Psi^{-1}(p_0) \in \mathbb{X}$.
		\For{$k=0,1,2,...,k_{\max}$}
		\State Solve the simplified Newton equation \eqref{eq:simplified_newton_equation},
		\begin{align*}
			v_k &= - (\tD \Psi(x_k)^\top W \tD \Psi(x_k))^\dag \tD \Psi(x_k)^\top W(p-z).
		\end{align*}
		\State Set
		\begin{align*}
			x_{k+1} &= x_k + v_k, \\
			\hat{p} &= \Psi(x_{k+1}).
		\end{align*}
		\If{$\vert v_k \vert \leq \text{threshold}$}
		\State \textbf{break}
		\EndIf
		\EndFor
	\end{algorithmic}
\end{algorithm}

Algorithm \ref{alg:newton_method} features guaranteed local convergence.
Specifically, if $x_* = \Psi^{-1}(p_*)$ is a critical point of $\vert p - z \vert^2_W$, then there is a neighbourhood $\mathcal{U}_{x_*}$ of $x_*$ such that, if $x_0 \in \mathcal{U}_{x_*}$ then the sequence $x_k$ generated by the algorithm converges at least quadratically to $x_*$ \cite[Theorem 6.3.2]{AbsiMaho09}.

The proposed algorithm requires an initial guess $x_0$, which can be obtained through two simple approaches in practice.
The first is to use a naive projection method such as coordinate projection.
That is, if the coordinate functions are included in the dictionary $\mathbb{V}$, so that $\Psi(x) = (x, \bar{\Psi}(x))$, where $\bar{\Psi} : \mathbb{X} \to \R^{M-d}$, then the initial guess can be chosen as the first $d$ coordinates of $p$, $x_0 := p_{1:d}$.
Another possibility, which does not rely on a particular structure of the dictionary, is to use the previous value of $x^*(t_k)$ to initialize the guess for the new value $x_0(t_{k+1}) = x^*(t_k)$.

\begin{figure}[htb]
	\centering
	\includegraphics[width=0.49\linewidth]{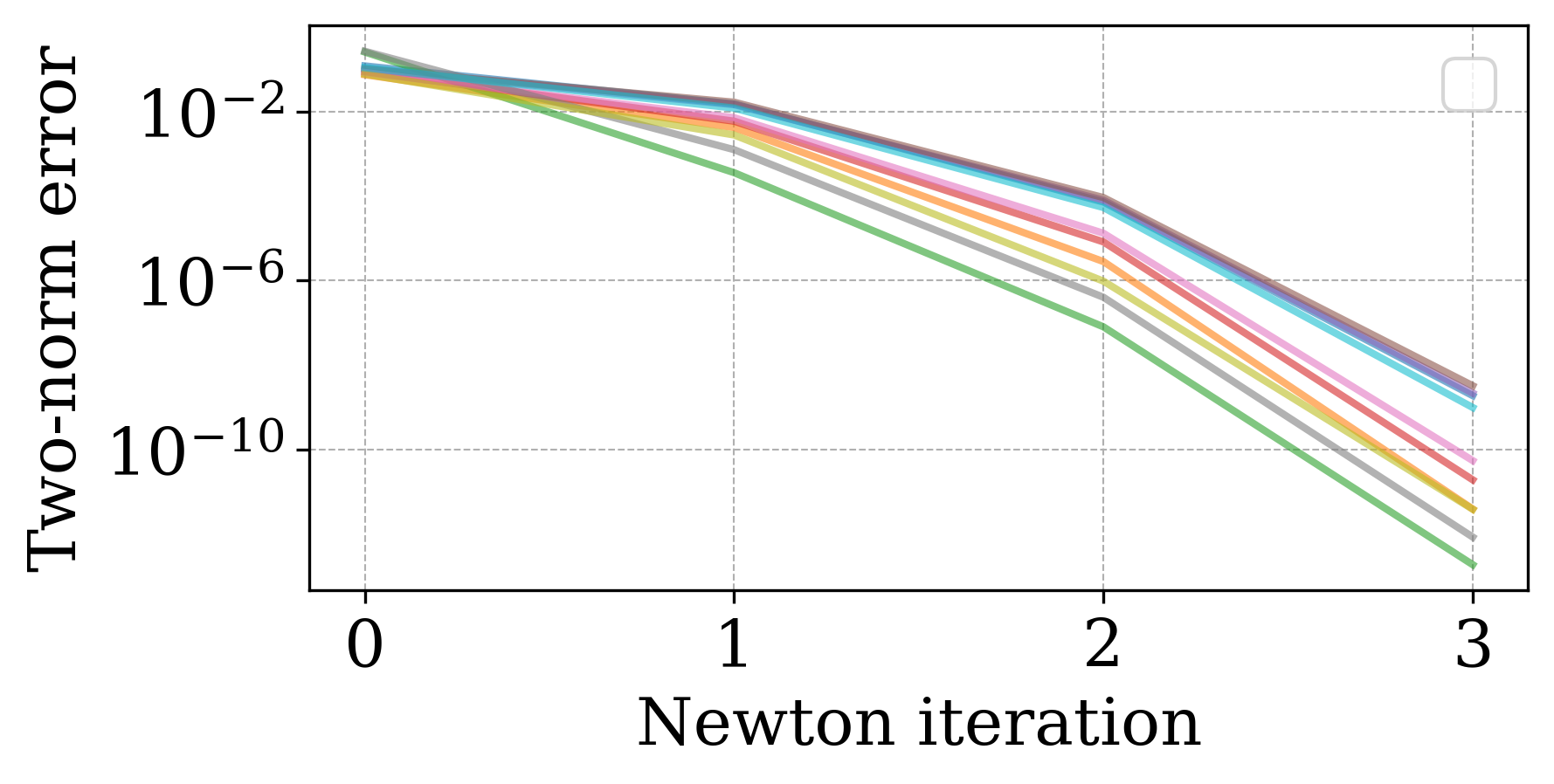}
	\includegraphics[width=0.49\linewidth]{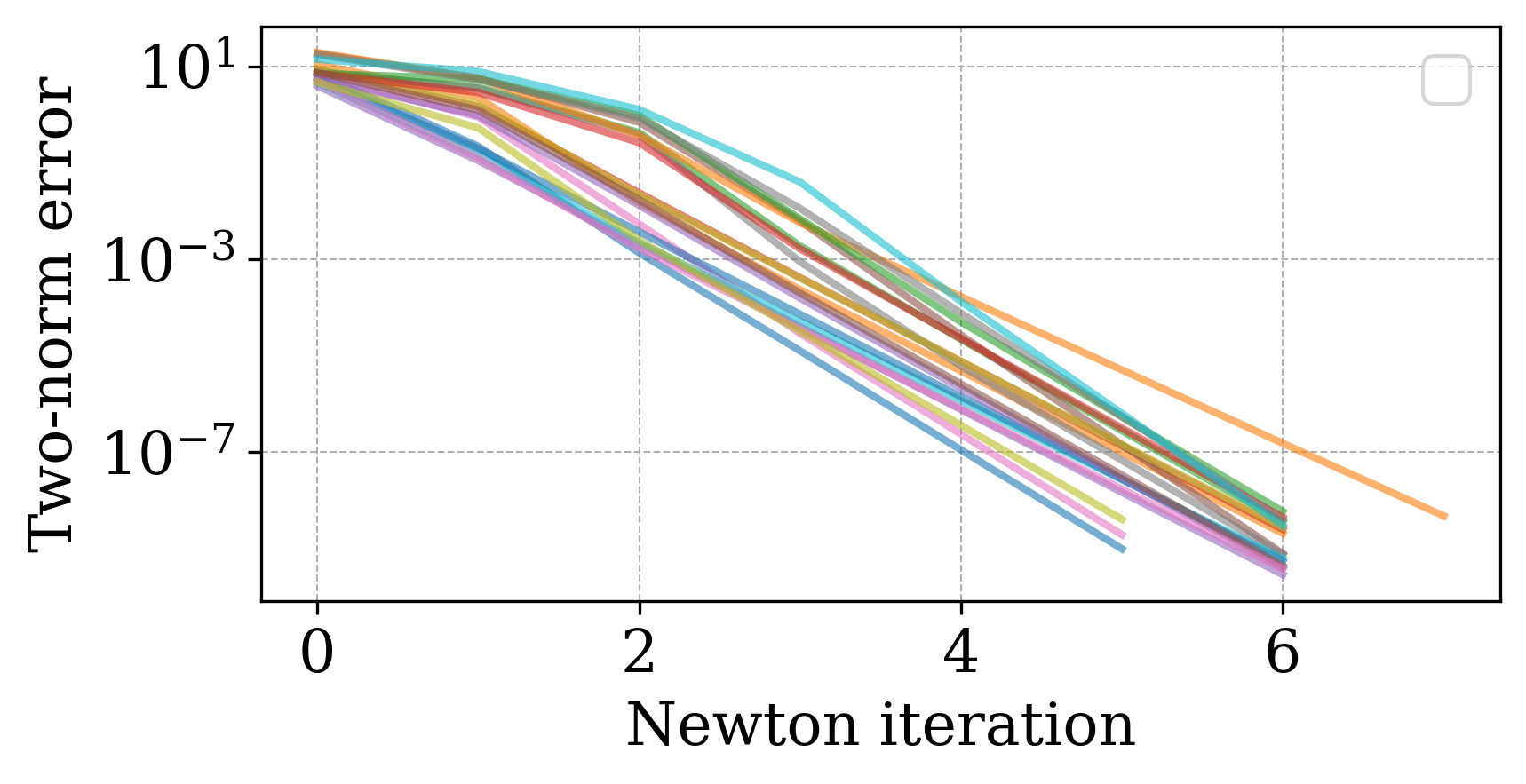}
	\caption{Error of Newton iteration for every 20th prediction step with time step $t=0.1$ for the Duffing oscillator (left) and the Lorenz example (right).}
	\label{fig:newton}
\end{figure}

In Figure~\ref{fig:newton}, we illustrate the superlinear convergence of the Riemannian Newton method for the projection. We observe that the tolerance of $10^{-8}$ is met already after a few time steps when using the solution of the previous time step as an initial guess.
We conclude this section by means of a remark considering globalization.
\begin{remark}
	The Newton method proposed in Algorithm \ref{alg:newton_method} can be extended to a trust region method such as \cite[Algorithm 10]{AbsiMaho09} to provide global convergence under some mild assumptions.
\end{remark}

\section{Conclusion}
\noindent In this paper, we proposed a reprojection method for Koopman-based predictions of nonlinear parameter-affine systems. Considering closest-point projections parametrized by weighting matrices mapping onto the consistent manifold of the compression, we have analyzed two particular choices corresponding to coordinate and maximum-likelihood projections and provided details for efficient implementation using an adaptive multistep scheme or a Riemannian Newton method. In various numerical examples, we have showcased the performance and increased robustness of Koopman-based predictions using these reprojections.

\section*{Declaration of the use of generative AI in the writing process}
\noindent During the preparation of this work, the authors used a
large language model (Claude; Anthropic) to assist with language editing and to draft a preliminary
formulation of Lemma~\ref{lem:centered}. All mathematical statements,
proofs and numerical results were developed, checked and verified by the
authors.

\bibliographystyle{abbrv}
\bibliography{references}

\end{document}